\documentclass[MNA]{mna} %

\usepackage{bbm}

\usepackage{subcaption}

\usepackage{pgfplots}
\pgfplotsset{width=10cm,compat=1.9}

\usepackage{environ}
\makeatletter
\newsavebox{\measure@tikzpicture}
\NewEnviron{scaletikzpicturetowidth}[1]{%
  \def\tikz@width{#1}%
  \def\tikzscale{1}\begin{lrbox}{\measure@tikzpicture}%
  \BODY
  \end{lrbox}%
  \pgfmathparse{#1/\wd\measure@tikzpicture}%
  \edef\tikzscale{\pgfmathresult}%
  \BODY
}

\SHORTTITLE{On a mean-field model of spiking neurons}

\TITLE{A mean-field model of Integrate-and-Fire neurons: non-linear stability of the stationary solutions
} 

\AUTHORS{%
  Quentin~Cormier\footnote{Inria, CMAP, CNRS, \'Ecole polytechnique, Institut Polytechnique de Paris, 91120 Palaiseau, France.
    \EMAIL{quentin.cormier@inria.fr}}
  }

\SHORTAUTHOR{Quentin Cormier}


\AMSSUBJ{60H10} 
\AMSSUBJSECONDARY{60K35; 45D05; 37A30} 

\SUBMITTED{November 22, 2023} 
\ACCEPTED{June 25, 2024} 

\VOLUME{4}
\YEAR{2024}
\PAPERNUM{1}
\DOI{10.46298/mna.12583}

\ABSTRACT{
We investigate a stochastic network composed of Integrate-and-Fire spiking neurons, focusing on its mean-field asymptotics. We consider an invariant probability measure of the McKean-Vlasov equation and establish an explicit sufficient condition to ensure the local stability of this invariant distribution.
Furthermore, we prove a conjecture proposed initially by J. Touboul and P. Robert regarding the bistable nature of a specific instance of this neuronal model. 
}

\newcommand*\indic[1]{\mathbbm{1}_{\{ #1 \}}}
\newcommand*\indica[1]{\mathbbm{1}_{ #1 }}

\def\cL{{\mathcal L}}
\def\cM{{\mathcal M}}

\def\cP{{\mathcal P}}

\def\C{\mathbb{C}}

\def\E{\mathbb{E}}

\def\R{\mathbb{R}}

\def\d{\mathrm{d}}

\def\Law{{\textstyle Law}}
\def\Lip{{\textstyle Lip}}

\providecommand{\abs}[1]{\lvert#1\rvert}
\providecommand{\norm}[1]{\lVert#1\rVert}

\begin{document}
	
\section{Introduction}
We consider a network of $N$ spiking neurons. Each neuron is characterized by its membrane potential $(X^{i, N}_t)_{t \geq 0}$. Each neuron emits ``spikes'' randomly, at a rate $f(X^{i, N}_t)$, which only depends on its membrane potential. The function $f: \mathbb{R} \rightarrow \mathbb{R}_+$ is deterministic. When a neuron spikes (say neuron $i$ spikes at time $\tau$), its potential is instantaneously reset to zero (we say zero is the \textit{resting value}) while the other neurons receive a small kick: 
\[  X^{i, N}_{\tau_+} = 0, \quad \text{ and } \quad \forall j \neq i, \quad  X^{j, N}_{\tau_+} = X^{j, 
N}_{\tau_-} + J^N_{i \rightarrow j}. \]
In this equation, the \textit{synaptic weight} $J^N_{i \rightarrow j}$ is a deterministic constant that models the interaction between the neurons $i$ and $j$. Finally, between the spikes, each neuron follows its own dynamics given by the scalar ODE:
\[  \dot{X}^{i, N}_t = b(X^{i, N}_t ), \]
where $b: \mathbb{R} \rightarrow \mathbb{R}$ is a deterministic function. We say that $b$ models the \textit{sub-threshold} dynamics of the neuron. We are interested in the dynamics of one particle (say $(X^{1, N}_t)$) in the limit where the number of particles $N$ goes to infinity. To simplify, we assume that the neurons are all-to-all connected with the same weight:
\[  \forall i,j, ~i \neq j \quad J^N_{i \rightarrow j} = \frac{J}{N}. \]
In this work, the deterministic constant $J$ is non-negative (when $f$ is non-decreasing, we say it is an \textit{excitatory network}).
At the initial time, all the neurons start with \textit{i.i.d.} initial conditions with law $\nu 
\in \mathcal{P}(\mathbb{R}_+)$. We assume that $b(0) \geq 0$, so the membrane 
potentials stay on $\mathbb{R}_+$. 

One expects \textit{propagation of chaos} to hold: as $N$ goes to infinity, any pair of neurons 
of the network (say $X^{1, N}_t$ and $X^{2, N}_t $) become more and more independent, and each 
neuron (say $(X^{1, N}_t)$) converges in law to the solution of the following McKean-Vlasov equation:
\begin{equation}
\label{eq:McKeanVlasov}
X^\nu_t = X^\nu_0 + \int_0^t{ b(X^\nu_s) \d s} + J \int_0^t{ \E f(X^\nu_s) \d s} - \int_0^t{ \int_{\R_+} X^\nu_{s-} \indic{z \leq f(X^\nu_{s-})} N(\d s, \d z) }.
\end{equation}
In this equation, $N(\d s, \d z)$ is a Poisson measure on $\mathbb{R}^2_+$ with intensity the Lebesgue measure $\d s \d z$. In addition, the initial condition $X^\nu_0$ has law $\nu \in \cP(\R_+)$  and is independent of the Poisson measure.
Informally, eq.~\eqref{eq:McKeanVlasov} can be understood in the following way: between the jumps, 
$(X^\nu_t)$ solves the ODE  $\dot{X}^\nu_t = b(X^\nu_t) + J \E f(X^\nu_t)$ and  $(X^\nu_t)$ jumps to zero at a rate 
$f(X^\nu_t)$.

This model of neurons is sometimes known in the literature as the ``Escape noise'', ``Noisy output'', ``Hazard rate'' model. We refer to \cite[Ch. 9]{GerstnetNeuronalDynamic} for a review. From a mathematical point of view, it was first introduced by \cite{de_masi_hydrodynamic_2015}, where it is described as a time-continuous version of the Galves–Löcherbach  model \cite{MR3055382}. Under few assumptions on $b$, $f$, and on the initial condition $\nu$, it is known that \eqref{eq:McKeanVlasov} is well-posed: see \cite{de_masi_hydrodynamic_2015} (with the assumption that the initial condition $\nu$ is compactly supported), \cite{fournier_toy_2016} (assuming only that $\nu$ has a first moment) and~\cite{CTV} (where a different proof is given, based on the renewal structure of the equation, see Theorem~\ref{th:stab} below). The convergence of the finite particle system $(X^{i, N}_t)$ to the solution of \eqref{eq:McKeanVlasov} is studied in \cite{fournier_toy_2016} where the rate of convergence, of the order $C(t) / \sqrt{N}$, is also given.

When the number of neurons is finite, the network $(X^{i, N}_t)$ is a Markov process (it is a Piecewise 
Deterministic Markov Process, see \cite{DavisPDMP}). So under quite general assumptions on $b$, 
 and $f$, this $\mathbb{R}^N_+$-valued process has a unique invariant probability measure, which is globally 
attractive. We refer to \cite{duarte_model_2014}, \cite{zbMATH06706511}, 
\cite{hodara_non-parametric_2018} and \cite{math7060518} for studies about the long time behavior 
of the finite particle system. 

The long time behavior of the solution of the limit equation~\eqref{eq:McKeanVlasov} is more complex, 
essentially because this is a McKean-Vlasov equation, and so it is not Markovian. In particular, 
\eqref{eq:McKeanVlasov} may have multiple invariant probability measures (see \cite{zbMATH06706511, CTV} and Section~\ref{sec:example} below for explicit 
examples). Even when the invariant probability measure is unique, it is not necessarily attractive. In \cite{drogoul_hopf_2017, CTV2}, the authors show that a Hopf bifurcation might appear when the interaction parameter $J$ varies, 
leading to periodic solutions of \eqref{eq:McKeanVlasov}.
The specific case $b \equiv  0$ is studied in \cite{fournier_toy_2016}. It is proved that for  all $J > 0$, there are precisely two invariant probability measures: the Dirac $\delta_0$, which is unstable,  and a non-trivial one, which is globally attractive.  This case $b \equiv 0$ is also investigated in \cite{veltz2}, where the authors prove that the non-trivial invariant measure is locally attractive with an exponential rate of convergence. Both \cite{fournier_toy_2016} and \cite{veltz2} rely on the Fokker-Planck PDE satisfied by the density of the solution of \eqref{eq:McKeanVlasov}. Finally, in \cite{CTV}, general conditions are given on $b$ and $f$ such that the McKean-Vlasov equation \eqref{eq:McKeanVlasov} admits a globally attractive invariant measure, assuming that the interaction parameter $J$ is small enough. For such weak enough interactions, similar results have been obtained for variants of this model, such as the time-elapsed model (see \cite{MR3850019}) or the Integrate-and-fire with a fixed deterministic threshold (see \cite{MR3190511}, \cite{MR3403402} and \cite{dumont_mean-field_2018} for another ``Poissonian'' variant). Finally, in \cite{zbMATH07481304}, the authors consider the case $b(x) = -x$, $f(x) = k \min(x, \lambda)$ for some constants $k, \lambda > 0$. They obtain a bistability result using a coupling method and study the metastable behavior of the particle system.

Understanding the long time behavior of \eqref{eq:McKeanVlasov} for an arbitrary interaction parameter $J$ is a challenging open question. In this work, we address the question of the local stability of an invariant probability measure of \eqref{eq:McKeanVlasov}, without assumptions on the size of the interactions $J$.
We say that $\nu_\infty \in \cP(\R_+)$ is an invariant probability measure of \eqref{eq:McKeanVlasov} if for all $t \geq 0$, it holds that the law of $X^{\nu_\infty}_t$, solution of \eqref{eq:McKeanVlasov}, is equal to $\nu_\infty$. We are interested in the stability of such an invariant probability measure $\nu_\infty$. Our main contribution is to provide an explicit criterion to decide stability. If this criterion is satisfied, then for any $\nu \in \cP(\R_+)$ sufficiently closes to $\nu_\infty$, the law of $X^\nu_t$ converges to $\nu_\infty$, at an exponential rate.

Recently, in \cite{cormier2023stability}, this question has been addressed for McKean-Vlasov with smooth coefficients driven by Brownian motions. It is shown that the stability is governed by the location of the zeros of a particular analytic function associated to the dynamics. We follow and adapt to \eqref{eq:McKeanVlasov} the strategy introduced in \cite{cormier2023stability}. 
One key difficulty is to identify a distance on probability measures that is well-suited to the particular structure of \eqref{eq:McKeanVlasov}. While the Wasserstein $W_1$ distance is used in \cite{cormier2023stability} for McKean-Vlasov driven by Brownian motions, we use here instead the ``Bounded Lipschitz'' distance, see eq.~\eqref{eq:bounded-lip-distance} below.
Using this distance, we first prove a stability result of \eqref{eq:McKeanVlasov} on a finite time interval $[0, T]$ (Theorem~\ref{th:stab}). Using the same distance, we then study the long time behavior of an auxiliary Markov process, see Proposition~\ref{prop:estimate-linear-process}. Combining these two results, we derive our main result, Theorem~\ref{th:stability}, which provides an explicit criterion to decide if an invariant probability measure is stable and quantifies the convergence.
We then identify in Proposition~\ref{prop:structural-stability} a structural condition on the coefficients, namely $f + b' \geq 0$, such that eq.~\eqref{eq:McKeanVlasov} has a unique, locally stable, invariant probability measure. This generalizes the results of \cite{fournier_toy_2016} and \cite{veltz2}, valid for $b \equiv 0$.

Finally, we study an explicit example where bistability occurs. We study the case:
\[ f(x) = x^2 \quad \text{ and } \quad b(x) = -x, \quad \forall x \in \R_+. \]
This example is studied numerically in \cite{zbMATH06706511}. In this work, the authors conjecture that \eqref{eq:McKeanVlasov} exhibits a phase transition: there is a parameter $J_* > 0$ such that for all $J < J_*$, \eqref{eq:McKeanVlasov} has a unique invariant probability measure (the Dirac mass at zero) while for all $J > J_*$, \eqref{eq:McKeanVlasov} has three invariant probability measures (the Dirac mass at zero and two other non-trivial probability measures). We provide a proof of this conjecture; see Proposition~\ref{prop:toy-model-multi-stab}.
Provided that $J > J_*$, an exact analysis of the stability of the two non-trivial stationary distributions seems impossible. However, a local analysis near the bifurcation points $J = J_*$ is possible. Then, in view of Theorem~\ref{th:stability}, we conjecture that one is stable and the other is unstable.
As the Dirac probability measure is stable (see \cite{zbMATH06706511}), we deduce that bistability is the main paradigm for $J > J_*$.

\textbf{Main notations}. 
We write $\cP(\R_+)$ for the space of probability measures on $\R_+$. 
Given $\nu \in \cP(\R_+)$ and $g: \R_+ \rightarrow \R$ a test function, we write $\langle g, \nu \rangle = \int_{\R_+} g(x) \nu(\d x)$. For $Z$ a random variable on $\R_+$, we write $\Law(Z) \in \cP(\R_+)$ for its probability law.
We denote by $\Lip_1(\R_+)$ the space of globally Lipschitz functions from $\R_+$ to $\R$ with Lipschitz norm bounded by $1$.
We equip $\cP(\R_+)$ with the following distance:  for all $\nu, \mu \in \cP(\R_+)$,
\begin{equation}
	\label{eq:bounded-lip-distance}
\norm{\nu-\mu}_0 := \sup \left\{ \int_{\R_+} g \d (\nu-\mu); \quad g \in \Lip_1(\R_+), \quad \sup_{x \in \R_+} |g(x)| \leq 1 \right\}.  
\end{equation}

\section{Main results}
Let $N(\d s, \d z)$ be a Poisson measure on $\R_+ \times \R_+$ with intensity the Lebesgue measure $\d s  \d z$.
Given $f: \R_+ \rightarrow \R_+$, $b: \R_+ \rightarrow \R$ and $J \geq 0$, we consider the McKean-Vlasov SDE \eqref{eq:McKeanVlasov}, 
where at the initial time, $X^\nu_0$ is distributed according to law $\nu \in \cP(\R_+)$. 
\begin{assumption} \label{ass:f-b-1}
	We assume that $b, f \in C^1(\R_+; \R_+)$ with $\norm{f}_\infty + \norm{f'}_\infty < \infty$ and that $b$ is globally Lipschitz, with $b(0) \geq 0$.
\end{assumption}
The condition $b(0) \geq 0$ ensures that the dynamics stays on $\R_+$.
We first show the following existence and stability result on a finite horizon time $[0, T]$:
\begin{theorem}\label{th:stab}
	Under Assumption~\ref{ass:f-b-1}, the mean-field equation~\eqref{eq:McKeanVlasov} has a unique path-wise solution for all $\nu \in \cP(\R_+)$. In addition, for all $T > 0$, there exists a constant $C_T$ such that for all $t \in [0, T]$:
	\[ \forall \nu, \mu \in \cP(\R_+), \quad \norm{\Law(X^\nu_t) - \Law(X^\mu_t)}_0 \leq C_T \norm{\nu - \mu}_0.  \]
\end{theorem}
Our second result concerns the ergodic behavior of the following associated linear equation. Let $\alpha > 0$. Denote by $(Y^{\alpha, \nu}_t)$ the solution of 
\begin{equation}
	\label{eq:linear-process}
	Y^{\alpha, \nu}_t =  Y^{\alpha, \nu}_0 + \int_0^t{ b(Y^{\alpha, \nu}_s) \d s} + \alpha t - \int_0^t{ \int_{\R_+} Y^{\alpha, \nu}_{s-} \indic{z \leq f(Y^{\alpha, \nu}_{s-})} N(\d s, \d z) },
\end{equation}
with $\Law(Y^{\alpha, \nu}_0) = \nu$. That is, we have ``frozen'' the non-linear interactions $J \E f(X^\nu_t)$ of \eqref{eq:McKeanVlasov}, and we have replaced it by the constant drift $\alpha$.

We denote by $\varphi^\alpha_t(x)$ the unique solution of the ODE
\[ \frac{\d }{\d t} \varphi^\alpha_t(x) = b(\varphi^\alpha_t(x)) + \alpha, \quad \text{ with } \quad \varphi^\alpha_0(x) = x. \]
In what follows, we assume that
\begin{assumption} \label{ass:f-b-2}
	The triple $(b,f, \alpha)$ satisfies
	\begin{enumerate}
		\item The jump rate is asymptotically lower bounded: 
			\begin{equation} \label{eq:underline_lambda}
				\underline{\lambda} := \liminf_{t \rightarrow \infty} \inf_{x \geq 0} \frac{1}{t} \int_0^t f(\varphi^\alpha_s(x)) \d s > 0. 
			\end{equation}
	\item There exists a constant $C$ such that for all $t \geq 0$ and $x, y \in \R_+$:
		\[ \abs{\varphi^\alpha_t(x) - \varphi^\alpha_t(y)} \leq C \abs{x-y}. \]
\end{enumerate}
\end{assumption}
\begin{remark}
	The constants $\underline{\lambda}$ and $C$ are allowed to depend on $\alpha$.
	The first point is satisfied if $f(x) \geq f_{\min} > 0$. The second point is satisfied if $b(x) = b_0  - b_1 x$ where $b_0, b_1 \geq 0$. Indeed, in that case, the flow is:
		\begin{equation} \label{eq:flow-linear-b} \varphi^\alpha_t(x) = \begin{cases} \left( x - \frac{b_0 + \alpha}{b_1}\right) e^{-b_1 t} + \frac{b_0 + \alpha}{b_1} \quad \text{ if } b_1 > 0, \\ x + (b_0 + \alpha) t \quad \text{ if } b_1 = 0. \end{cases}
\end{equation}
\end{remark}
We have:
\begin{proposition}
	\label{prop:estimate-linear-process}
	Under Assumption~\ref{ass:f-b-1} and~\ref{ass:f-b-2}, there exists $C_* > 0$ and $\lambda_* \in (0, \underline{\lambda})$ such that for all $\nu, \mu \in \cP(\R_+)$, and for all $t \geq 0$
	\[ \norm{\Law(Y^{\alpha, \nu}_t) - \Law(Y^{\alpha, \mu}_t)}_0 \leq C_* e^{-\lambda_* t} \norm{\nu-\mu}_0. \]
\end{proposition}
Therefore, \eqref{eq:linear-process} has a unique invariant probability measure. As shown in \cite{CTV}, this invariant probability measure is:
\begin{equation}
	\label{eq:invariant-measure}
\nu^\alpha_\infty(x) \d x = \frac{\gamma(\alpha)}{b(x) + \alpha} \exp \left(- \int_0^x \frac{f(y)}{b(y) + \alpha} \d y \right) \indica{[0, \sigma_\alpha)}(x) \d x, 
\end{equation}
where the support $\sigma_\alpha \in (0, \infty]$ is given by
\[ \sigma_\alpha = \inf\{y \geq 0; \quad b(y) + \alpha = 0 \} = \lim_{t \rightarrow \infty} \varphi^\alpha_t(0). \]
The changes of variables $y = \varphi^\alpha_u(0)$ and $x = \varphi^\alpha_t(0)$ in \eqref{eq:invariant-measure} show that for any bounded measurable test function $g$,
\[ \langle g, \nu^\alpha_\infty \rangle = \gamma(\alpha) \int_0^\infty g(\varphi^\alpha_t(0)) \exp\left(-\int_0^t f(\varphi^\alpha_u(0)) \d u \right) \d t. \]
In particular,  choosing $g = f$, we find that $\gamma(\alpha)$ satisfies:
\[  \gamma(\alpha) = \int_{\R_+} f(x) \nu^\alpha_\infty(x) \d x.  \]
We say that $\nu_\infty$, an invariant probability measure of \eqref{eq:McKeanVlasov}, is non-trivial if $J \langle f, \nu_\infty\rangle \neq 0$.
Note that there is a one-to-one correspondence between the non-trivial invariant probability measures of \eqref{eq:McKeanVlasov} and the $\alpha > 0$ satisfying:
\[ \alpha = J \gamma(\alpha). \]
Let $\nu_\infty$ be a non-trivial invariant probability measure of \eqref{eq:McKeanVlasov}. Let $\alpha = J \langle f, \nu_\infty \rangle$ such that $\nu_\infty = \nu^\alpha_\infty$. We assume that the triple $(b, f, \alpha)$ satisfies Assumption~\ref{ass:f-b-2}. To state our main result on the stability of $\nu_\infty$, we need to introduce two notations.
We define for all $t \geq 0$
\begin{align}
H_\alpha(t) &:= \exp \left( - \int_0^t f(\varphi^\alpha_s(0)) \d s \right) \label{eq:def-of-Halpha}, \\
\Psi_\alpha(t) &:= \alpha \int_0^\infty H_\alpha(t+u) \frac{f(\varphi^\alpha_{t+u}(0)) - f(\varphi^\alpha_u(0)) }{b(\varphi^\alpha_u(0)) + \alpha} \d u. \label{eq:def-of-Psi-alpha}
\end{align}
Let $\lambda_* > 0$ be the constant given by Proposition~\ref{prop:estimate-linear-process}.
\begin{assumption}
	\label{ass:final}
	Assume that $f \in C^2(\R_+; \R_+)$ with $\norm{f^{(k)}}_\infty < \infty$ for $k \in \{0, 1, 2 \}$.
	Assume $b \in C^2(\R_+; \R)$ with $b(0) \geq 0$ and $\norm{b'}_\infty + \norm{b''}_\infty < \infty$.
	Finally, assume there exists a constant $\lambda_\alpha \in (0, \lambda_*)$ such that 
	\[ \sup_{t \geq 0} e^{\lambda_\alpha t} \int_0^\infty H_\alpha(t+u) \left| \frac{f(\varphi^\alpha_{t+u}(0)) - f(\varphi^\alpha_u(0)) }{b(\varphi^\alpha_u(0)) + \alpha} \right| \d u< \infty. \]
\end{assumption}
\begin{remark}
	This last estimate holds if $b(x) = b_0 - b_1 x$ with $b_1 \geq 0$. This follows from the explicit expression of $\varphi^\alpha_t(0)$, see \eqref{eq:flow-linear-b}, and from the fact that $f$ is globally Lipschitz.
\end{remark}
We write $D_\alpha := \{z \in \C,~\Re(z) > -\lambda_\alpha\}$. We consider for $z \in D_\alpha$
\[ \hat{H}_\alpha(z) := \int_0^\infty e^{-zt} H_\alpha(t) \d t,\] 
the Laplace transform of $H_\alpha$. Similarly, let $\hat{\Psi}_\alpha(z)$ be the Laplace transform of $\Psi_\alpha$. By assumption, the two functions are analytic on $D_\alpha$. Our main result is
\begin{theorem}
	\label{th:stability}
	Let $\nu_\infty$ be a non-trivial invariant probability measure of \eqref{eq:McKeanVlasov}. Let $\alpha = J \langle f, \nu_\infty \rangle \in \R^*_+$. In addition to Assumptions~\ref{ass:f-b-2} and~\ref{ass:final}, assume that
	\begin{equation} \label{eq:spectral-assumption}
	-\lambda'_\alpha := \sup \{\Re(z); \quad z \in D_\alpha, \quad \hat{H}_\alpha(z) = \hat{\Psi}_\alpha(z)  \} < 0. 
\end{equation}
Then there exists $C, \epsilon > 0$ and $\lambda \in (0, \lambda'_\alpha)$ such that for all initial condition $\nu \in \cP(\R_+)$ with $\norm{\nu - \nu_\infty}_0 < \epsilon$, it holds that
	\[ \forall t \geq 0, \quad \norm{\Law(X^\nu_t) - \nu_\infty}_0 \leq C e^{-\lambda t} \norm{\nu-\nu_\infty}_0.  \] 
\end{theorem}
We provide in Section~\ref{sec:example} an explicit example with multiple invariant probability measures.
The spectral assumption $\lambda'_\alpha > 0$ is automatically satisfied if the structural condition $f + b' \geq 0$ holds. More precisely,
\begin{proposition}
	\label{prop:structural-stability}
	We have:
	\begin{enumerate}
		\item Let $\nu_\infty$ be a non-trivial invariant probability measure of \eqref{eq:McKeanVlasov}. Let $\alpha = J \langle f, \nu_\infty \rangle \in \R^*_+$. In addition to Assumption~\ref{ass:f-b-2} and~\ref{ass:final}, assume that $f + b' \geq 0$ on $[0, \sigma_\alpha)$. Then $\lambda'_\alpha > 0$, and so $\nu_\infty$ is locally stable.
		\item Grant Assumption~\ref{ass:f-b-1}. Assume moreover that $b \in C^1(\R_+)$, $b(0) > 0$ and that $f + b' \geq 0$ on $\R_+$. Then \eqref{eq:McKeanVlasov} has exactly one non-trivial invariant probability measure.
	\end{enumerate}
\end{proposition}
\section{Proofs}
\subsection{Notations}
Let $T > 0$. 
Given $a \in C([0, T]; \R_+)$, we denote by $Y^{a, \nu}_{t,s}$ the solution of the linear non-homogeneous SDE
\begin{equation}\label{eq:linear-non-homogeneous}
	Y^{a, \nu}_{t,s} = Y^{a, \nu}_{s,s} + \int_s^t{ b(Y^{a, \nu}_{u, s}) \d u} +  \int_s^t{ a_u \d u} - \int_s^t{ \int_{\R_+} Y^{a, \nu}_{u-, s} \indic{z \leq f(Y^{a, \nu}_{u-, s})} N(\d u, \d z) },
\end{equation}
where at time $s$, $\Law(Y^{a, \nu}_{s,s}) = \nu$.
We let $\varphi^{a}_{t, s}(x)$ be the solution of the ODE
\[ \frac{\d}{\d t} \varphi^a_{t, s}(x) = b(\varphi^a_{t,s}(x)) + a_t, \quad \varphi^a_{s, s}(x) = x. \]
As in \cite{CTV}, we denote by $K_a^\nu(t,s)$ the density of the first jump of $(Y^{a, \nu}_{t, s})_{t \geq s}$:
\[ K^\nu_a(t, s) := \int_{\R+}f(\varphi^a_{t,s}(x)) \exp \left(-\int_s^t{ f(\varphi^a_{u, s}(x)) \d u} \right) \nu(\d x).   \]
Similarly, let $H^\nu_a(t, s)$ be the survival function of the first jump:
\[  H^\nu_a(t, s) := \int_{\R+} \exp \left(-\int_s^t{ f(\varphi^a_{u, s}(x)) \d u} \right) \nu(\d x).   \]
When $a$ does not depend on $t$, that is $a \equiv \alpha \in \R_+$, we write for all $t, x \geq 0$:
\begin{equation}
	\label{eq:H-x-alpha-and-K-x-alpha}
	H^x_\alpha(t) = H^{\delta_x}_\alpha(t, 0), \quad K^x_\alpha(t) = K^{\delta_x}_\alpha(t, 0), \quad \varphi^\alpha_t(x) = \varphi^\alpha_{t, 0}(x). 
\end{equation}
\subsection{Proof of Theorem~\ref{th:stab}}
\begin{lemma}\label{lem:estimate-Ka-Kb-gen}
	There exists a constant $C_T$ such that for all $g \in C^1(\R_+)$, for all $a,\tilde{a} \in C([0, T]; \R_+)$, for all $0 \leq s \leq t \leq T$,
\begin{align*} &\left|  \int_{\R_+} g(\varphi^{a}_{t, s}(x)) H^{\delta_x}_a(t, s) \nu(\d x) -  \int_{\R_+} g(\varphi^{\tilde{a}}_{t, s}(x)) H^{\delta_x}_{\tilde{a}}(t, s) \mu(\d x) \right| \\ & \quad \quad  \leq C_T( \norm{g}_\infty + \norm{g'}_\infty ) \left[  \int_s^t{ \abs{a_u -  \tilde{a}_u} \d u} +  \norm{\nu-\mu}_0 \right]. \end{align*}
\end{lemma}
\begin{proof}
	Write $F_a(x) = g(\varphi^a_{t, s}(x)) H^{\delta_x}_a(t, s)$. Using the explicit formula satisfied by the survival function $H^{\delta_x}_a(t, s)$, we find that 
	$F_a$ is $C^1(\R_+)$ with $\norm{F_a}_\infty \leq \norm{g}_\infty$ and 
	\[ \norm{F'_a}_\infty \leq \left[ \norm{g'}_\infty + \norm{g}_\infty \norm{f'}_\infty \right] e^{T \norm{b'}_\infty} . \]
	Therefore, for $C_T =  (1+ e^{T \norm{b'}_\infty}) (1+ \norm{f'}_\infty)$, we have:
	\[ \abs{\langle F_a, \nu-\mu \rangle} \leq C_T (\norm{g}_\infty + \norm{g'}_\infty) \norm{\nu - \mu}_0. \]
	In addition, using Gronwall's inequality, we find that
	\[ \abs{\varphi^a_{t, s}(x) - \varphi^{\tilde{a}}_{t, s}(x)} \leq e^{T \norm{b'}_\infty} \int_s^t \abs{a_u - \tilde{a}_u} \d u. \]
	Using the explicit formula of $F_a$, we deduce that there exists another constant $C_T$ such that
	\[ \abs{F_a(x) - F_{\tilde{a}}(x)} \leq C_T (\norm{g}_\infty + \norm{g'}_\infty) \int_s^t \abs{a_u - \tilde{a}_u} \d u. \]
	Altogether, we deduce the result.
\end{proof}
Let $r^\nu_a(t,s) := \E f(Y^{a, \nu}_{t,s})$. It holds that the jump rate $r^\nu_a(t,s)$ and $K^\nu_a(t, s)$ are linked by the following Volterra integral equation \cite{CTV}:
\begin{lemma} \label{lem:volterra-eq} 
	It holds that:
	\[ r^\nu_a(t,s) = K^\nu_a(t, s) + \int_s^t  r^{\delta_0}_a(t, u) K^\nu_a(u, s) \d u. \]
\end{lemma}
\begin{proof}
 Let $t \geq s$ and $\tau_s := \inf\{ u > s, Y^{a, \nu}_{u-, s} \neq Y^{a, \nu}_{u, s} \}$ be the time of the first jump of $Y^{a, \nu}_{\cdot, s}$ after $s$. The law of $\tau_s$ is $K^\nu_{a}(u, s) \d u$. 
We have
\[ r^\nu_a(t, s) =  \E f(Y^{a, \nu}_{t, s}) =  \E f(Y^{a, \nu}_{t, s}) \indic{\tau_s \geq t} +  \E f(Y^{a, \nu}_{t, s}) \indic{\tau_s  \in (s, t)}. \]
For a fixed initial condition $x$, it holds that $Y^{a, \delta_x}_{t,s} = \varphi^{a}_{t, s}(x)$ under the event $\{\tau_s > t\}$. Therefore,  the first term is equal to
\[  \int_{\R_+} f(\varphi^a_{t, s}(x)) H^{\delta_x}_a(t, s) \nu(\d x) = K^\nu_a(t, s). \]
Using the strong Markov property at time $\tau_s$ and using that the process is reset to $0$ after this jump, we find that the second term is equal to
\[   \E f(Y^{a, \nu}_{t, s}) \indic{\tau_s  \in (s, t)} = \E f(Y^{a, \delta_0}_{t, \tau_s} \indic{\tau_s \in (s, t)})=  \int_s^t{ r_{a}^{\delta_0}(t, u) K^\nu_{a}(u, s) \d u}. \]
Altogether, we deduce the result.
\end{proof}
More generally, by the exact same argument, we have
\begin{lemma}
	\label{lem:flow-SDE-volterra}
	Let $g: \R_+ \rightarrow \R$ be a bounded test function. It holds that for all $t \geq s$, 
	\[ \E g(Y^{a, \nu}_{t, s}) = \int_{\R_+}g(\varphi^{a}_{t, s}(x)) H^{\delta_x}_a(t, s) \nu(\d x) + \int_0^t K^\nu_a(u, s) \E g(Y^{a, \delta_0}_{t, u}) \d u.  \]
\end{lemma}
Exploiting the Volterra integral equation of Lemma~\ref{lem:volterra-eq}, we deduce that:
\begin{lemma}
	\label{lem:stab-non-homogeneous}
	There exists a constant $C_T$ such that for all $a, \tilde{a} \in C([0, T]; \R_+)$, it holds that
	\[   \forall s \leq t \leq T, \quad \abs{r^\nu_a - r^\mu_{\tilde{a}}}(t, s)  \leq C_T \int_s^t{ \abs{a_u - \tilde{a}_u} \d u} + C_T \norm{\nu-\mu}_0. \] 
\end{lemma}
\begin{proof}
	We first prove the inequality when $\nu = \mu = \delta_0$.  
	To simplify the notations, we also denote by $r_a(t, s) := r^{\delta_0}_a(t, s)$ and $K_a(t, s) := K^{\delta_0}_a(t, s)$.
Let $\Delta(t,s) := r_a(t, s) - r_{\tilde{a}}(t,s)$. Using Lemma~\ref{lem:volterra-eq}, we have
	\[ \Delta(t,s) = (K_a - K_{\tilde{a}})(t, s) + \int_s^t{ (K_a - K_{\tilde{a}})(u,s) r_a(t, u) \d u } + \int_s^t{ K_{\tilde{a}}(u,s) \Delta(t,u) \d u}. \]
	As $|r_a(t, s)| \leq \norm{f}_\infty$ and $|K_a(t, s)| \leq \norm{f}_\infty$, this shows that $(t, s) \mapsto \Delta(t, s)$ is continuous. In addition we have:
	\[ 	\abs{\Delta (t,s)} \leq C_T \int_s^t{ \abs{a_u - \tilde{a}_u} \d u }  + C_T \int_s^t{ \abs{\Delta(t, u)} \d u}.  \]
	We conclude by using the Grönwall's inequality. The extension to arbitrary $\nu, \mu \in \cP(\R_+)$ follows from Lemma~\ref{lem:estimate-Ka-Kb-gen} (with $g = f$) and Lemma~\ref{lem:flow-SDE-volterra}.
\end{proof}
By similar arguments, we have, using Lemma~\ref{lem:estimate-Ka-Kb-gen} and Lemma~\ref{lem:flow-SDE-volterra}:
	\begin{lemma}
		\label{lem:stability-Y-g}
		There exists a constant $C_T$ such that for all $a, \tilde{a} \in C([0, T]; \R_+)$, for all $g \in C^1(\R_+)$:	
		\[   \forall s \leq t \leq T, \quad  \abs{ \E g(Y^{a, \nu}_{t, s}) - \E g(Y^{\tilde{a}, \mu}_{t, s})} \leq C_T (\norm{g}_\infty + \norm{g'}_\infty) \left[   \int_s^t{ \abs{a_u - \tilde{a}_u} \d u} + \norm{\nu-\mu}_0 \right]. \] 
	\end{lemma}
We now give the proof of Theorem~\ref{th:stab}. 
Existence is proven exactly as in \cite{CTV}. We prove the stated stability estimate, which implies uniqueness.
 Let $(X^\nu_t)$ and $(X^\mu_t)$ be two solutions of \eqref{eq:McKeanVlasov} starting from laws $\nu$ and $\mu$.
Let $a_t = J \E f(X^\nu_t)$ and $\tilde{a}_t = J \E f(X^\mu_t)$. Then,  $a \in C([0, T]; \R_+)$. Indeed, by Ito's formula, we have
\[ a_t = J \E f(X^\nu_t) = J \E f(X^\nu_0) + \int_0^t{ J \E f'(X^\nu_s) (b(X^\nu_s) + a_s) \d s} + \int_0^t J \E (f(0) - f(X^\nu_s)) f(X^\nu_s) \d s. \]
In addition, $(X^\nu_t)$ is a solution of \eqref{eq:linear-non-homogeneous} with $a$. The same holds for $(X^\mu_t)$ with $\tilde{a}$. Therefore, by Lemma~\ref{lem:stab-non-homogeneous} (with $s = 0$), we deduce that
\[ |a_t - \tilde{a}_t| \leq C_T \int_0^t{ |a_u - \tilde{a}_u| \d u} + C_T \norm{\nu-\mu}_0. \]
By Grönwall's inequality, we deduce that
\[ \sup_{t \in [0, T]}|a_t - \tilde{a}_t| \leq C_T e^{C_T} \norm{\nu - \mu}_0. \]
The stability estimate of Theorem~\ref{th:stab} then follows from Lemma~\ref{lem:stability-Y-g}. \qed
\subsection{Proof of Proposition~\ref{prop:estimate-linear-process}}
Recall that the Markov process $(Y^{\alpha, \nu}_t)$ is defined by \eqref{eq:linear-process}. 
We also use the notation $\E_x g(Y^\alpha_t) := \E g(Y^{\alpha, \delta_x}_t)$, for all $x \geq 0$.
The first step is to prove that $r_\alpha(t) := \E_0 f(Y^{\alpha}_t)$ converges to $\gamma(\alpha)$ at an exponential speed:
\begin{lemma}
	Under assumptions~\ref{ass:f-b-1} and~\ref{ass:f-b-2}, there is a constant $\theta_\alpha > 0$ such that
	\[  \sup_{t \geq 0} \abs{r_\alpha(t) - \gamma(\alpha)} e^{\theta_\alpha t}  < \infty. \]
\end{lemma}
\begin{proof}
	By Lemma~\ref{lem:volterra-eq}, $r_\alpha$ is the solution of a Volterra convolution equation.
	Therefore, the strategy is to use the Laplace transform to deduce the asymptotic behavior of $t \mapsto r_\alpha(t)$ from the location of the poles of $\widehat{r}_\alpha(z)$. The arguments can be found in~\cite{CTV}. We only use here that $f$ is $C^1$, $f' b$ and $f^2$ are bounded and that
	$\liminf_{t \rightarrow \infty} \frac{1}{t} \int_0^t f(\varphi^\alpha_s(0)) \d s > 0.$
\end{proof}
We then show that we have convergence in total variation norm:
\begin{lemma}
	There are constants $C, \theta_\alpha > 0$ such that for all $g \in C(\R_+)$ with $\norm{g}_\infty < \infty$: 
\[ \sup_{x \geq 0}\abs{\E_x g(Y^\alpha_t) - \langle g, \nu^\alpha_\infty \rangle} \leq C \norm{g}_\infty e^{-\theta_\alpha t}, \quad \forall t \geq 0. \]
\end{lemma}
\begin{proof}
	Recall that $H_\alpha(t)$ is defined by \eqref{eq:def-of-Halpha} and that $K_\alpha(t) := -\frac{\d }{\d t} H_\alpha(t)$.
	We first show the result for $x = 0$. By Lemma~\ref{lem:flow-SDE-volterra},
	\[ \E_0 g(Y^\alpha_t) = g(\varphi^\alpha_t(0)) H_\alpha(t) + \int_0^t{ K_\alpha(u) \E_0 g(Y^\alpha_{t-u}) \d u}.  \]
	We solve this Volterra equation and find:
	\[ \E_0 g(Y^\alpha_t) = g(\varphi^\alpha_t(0)) H_\alpha(t) + \int_0^t r_\alpha(t-u) g(\varphi^\alpha_u(0)) H_\alpha(u) \d u.     \]
	We used here that $r_\alpha$ is the resolvent of $K_\alpha$, see \cite[Ch. 2]{gripenberg_volterra_1990}.
	We write $r_\alpha(t) = \gamma(\alpha) + \xi_\alpha(t)$ with $\abs{\xi_\alpha(t)} \leq Ce^{-\theta_\alpha t}$. We deduce that
	\begin{align*} \E_0 g(Y^\alpha_t) - \gamma(\alpha) \int_0^\infty  g(\varphi^\alpha_u(0)) H_\alpha(u) \d u &= g(\varphi^\alpha_0(t)) H_\alpha(t) + \int_0^t{ \xi_\alpha(t-u) g(\varphi^\alpha_0(u)) H_\alpha(u) \d u } \\
		& \quad - \gamma(\alpha) \int_t^\infty g(\varphi^\alpha_0(u)) H_\alpha(u)  \d u =: A_1+A_2+A_3.
	\end{align*}
	Recall that $\nu^\alpha_\infty$ is given by \eqref{eq:invariant-measure}. The change of variable $x = \varphi^\alpha_u(0)$ shows that
	\[   \gamma(\alpha) \int_0^\infty  g(\varphi^\alpha_u(0)) H_\alpha(u) \d u = \langle g, \nu^\alpha_\infty \rangle.   \]
	Using \eqref{eq:underline_lambda}, we deduce that there exists constants $C, \lambda > 0$ such that
	\[  \forall t \geq 0, \quad H_\alpha(t) \leq C e^{-\lambda t}. \]
	Without loss of generality, we can choose $\theta_\alpha < \lambda/2$. 
	Therefore, $|A_1| \leq C \norm{g}_\infty e^{-2 \theta_\alpha t}$. Similarly, 
	\[ |A_2| \leq C^2 \norm{g}_\infty \int_0^t e^{-\theta_\alpha (t-u) }e^{-2 \theta_\alpha u} \d u \leq C^2 \norm{g}_\infty e^{-\theta_\alpha t} \int_0^\infty e^{-\theta_\alpha u} \d u. \]
	Finally, 
	\[ |A_3| \leq C \norm{f}_\infty \norm{g}_\infty e^{-\theta_\alpha t} \int_0^\infty {e^{-\theta_\alpha u} \d u}.  \]
Altogether, there exists another constant $C$ such that:
	\[ \abs{\E_0 g(Y^\alpha_t) - \langle g, \nu^\alpha_\infty \rangle} \leq C \norm{g}_\infty e^{-\theta_{\alpha} t}.  \]
	Finally, we treat the general case. 
	Recall that $H^x_\alpha(t)$ and $K^x_\alpha(t)$ are defined by \eqref{eq:H-x-alpha-and-K-x-alpha}.
	For all $x \geq 0$, we have by Lemma~\ref{lem:flow-SDE-volterra}:
	\[ \E_x g(Y^\alpha_t)  = g(\varphi^\alpha_t(x)) H^x_\alpha(t) + \int_0^t{ K^x_\alpha(u) \E_0 g(Y^\alpha_{t-u})\d u}, \]
	and so, using that $\int_0^\infty K^x_\alpha(u) \d u =1$, it holds that
	\[  \E_x g(Y^\alpha_t) - \langle g, \nu^\alpha_\infty \rangle =  g(\varphi^\alpha_t(x)) H^x_\alpha(t) + \int_0^t{ K^x_\alpha(u) (\E_0 g(Y^\alpha_{t-u})  - \langle g, \nu^\alpha_\infty \rangle  )\d u} + \langle g, \nu^\alpha_\infty \rangle H^x_\alpha(t).    \]
	Therefore, the stated estimate is deduced from the case $x = 0$.
\end{proof}
Recall that the constant $\underline{\lambda} > 0$ is defined by \eqref{eq:underline_lambda}.
The next step is to prove that
\begin{lemma}
	There exists a constant $C$ such that for all $g \in C^1(\R_+)$ with $\norm{g}_\infty \leq 1$ and $\norm{g'}_\infty \leq 1$, for all $x, y \geq 0$:
	\[ \abs{g(\varphi^\alpha_t(x)) H^x_\alpha(t) - g(\varphi^\alpha_t(y)) H^y_\alpha(t)} \leq C \abs{x-y} e^{-(\underline{\lambda}/2) t}. \] 
\end{lemma}
\begin{proof}
We first show the result when $g \equiv 1$. We use the inequality $\abs{e^{-A} - e^{-B}} \leq e^{-\min(A, B)} \abs{A-B}$, valid for all $A, B \geq 0$. Using Assumption~\ref{ass:f-b-2}, we deduce that 
	\[ \abs{ H^x_\alpha(t) - H^y_\alpha(t) } \leq C e^{- (2\underline{\lambda}/3) t} \int_0^t \abs{f(\varphi^\alpha_u(x)) - f(\varphi^\alpha_u(y)) } \d u. \]
	Using that $f$ is globally Lipschitz and that $|\varphi^\alpha_t(x) - \varphi^\alpha_t(y)| \leq C |x-y|$, we deduce the stated inequality.
The general case is deduced similarly, as $g$ and $g'$ are assumed to be bounded by one. 
\end{proof}
Finally, we give the proof of Proposition~\ref{prop:estimate-linear-process}. 
In what follows, the constant $\lambda > 0$ might decrease from line to line.
Let $g \in C^1(\R_+)$ with $\norm{g}_\infty \leq 1$ and $\norm{g'}_\infty \leq 1$. We write
\begin{align*}
	\E_x g(Y^\alpha_t) - \E_y g(Y^\alpha_t) &= g(\varphi^\alpha_t(x)) H^x_\alpha(t) - g(\varphi^\alpha_t(y)) H^y_\alpha(t) \\
	& \qquad+ \int_0^t{ (K^x_\alpha(u) - K^y_\alpha(u)) \E_0 g(Y^\alpha_{t-u})  \d u}.
\end{align*}
The first term is bounded by $C |x-y| e^{-\underline{\lambda} t}$ by the previous lemma. In addition, 
\[ \forall t \geq u, \quad \abs{\E_0 g(Y^\alpha_{t-u}) - \langle g, \nu^\alpha_\infty \rangle} \leq C e^{-\lambda (t-u)}. \] 
So, using that $\int_0^\infty (K^x_\alpha(u) - K^y_\alpha(u)) \d u$ = 0, we find that the second term is bounded by
\[ \abs{H^x_\alpha(t) - H^y_\alpha(t)} \abs{ \langle g, \nu^\alpha_\infty \rangle} + C |x-y| e^{-\lambda t }.   \]
Altogether, we deduce that there is a constant $C > 0$ and $\lambda > 0$ such that for all $x, y$:
\[ \abs{ \E_x g(Y^\alpha_t)  - \E_y g(Y^\alpha_t) } \leq C |x-y| e^{-\lambda t}.   \]
We define: 
\[ v_t(x) := (\E_x g(Y^\alpha_t) - \langle g, \nu^\alpha_\infty \rangle ) e^{\lambda t}. \]
In view of Lemma~\ref{lem:flow-SDE-volterra}, it holds that $v_t \in C^1(\R_+)$. 
By the previous results, we have for some constant $C$ (independent of $g$):
\[ \norm{v_t}_\infty + \norm{v'_t}_\infty \leq C. \]
So
\[ \abs{ \langle v_t, \nu - \mu \rangle } \leq C \norm{\nu - \mu}_0.  \]
In other words, by the Markov property: 
\[ \norm{\Law(Y^{\alpha, \nu}_t) - \Law(Y^{\alpha, \mu}_t)}_0 \leq Ce^{-\lambda t} \norm{\nu-\mu}_0. \]
This ends the proof. \qed
\subsection{Reformulation of the spectral assumption}
In this section, we reformulate the spectral assumption~\eqref{eq:spectral-assumption}.
We recall that 
\begin{align*}
	H^y_\alpha(t) &=  \exp\left(-\int_0^t{f(\varphi^\alpha_u(y)) \d u}\right), \quad 
	H_\alpha(t) = H^0_\alpha(t), \\
	\Psi_\alpha(t) &:= \alpha \int_0^\infty H_\alpha(t+u) \frac{f(\varphi^\alpha_{t+u}(0)) - f(\varphi^\alpha_u(0))}{b(\varphi^\alpha_u(0)) + \alpha} \d u.
\end{align*}
The first step is to show that
\begin{lemma}
	It holds that
	\begin{equation} 
		\label{eq:formula-for-Psi}
	\Psi_\alpha(t) = -J \int_0^{\sigma_\alpha} \left[ \frac{\d }{\d y} H^y_\alpha(t) \right] \nu^\alpha_\infty(\d y). \end{equation}
\end{lemma}
\begin{proof}
	First, note that the function $y \mapsto \varphi^\alpha_u(y)$ is $C^1(\R_+)$ with
	\begin{equation} \label{eq:derivative-flow} \frac{\d}{\d y} \varphi^\alpha_u(y)  = \frac{b(\varphi^\alpha_u(y)) + \alpha}{b(y) + \alpha}. \end{equation}
	Indeed, both the left-hand-side and the right-and-side of \eqref{eq:derivative-flow} satisfies the ODE: $\partial_u \psi_u = b'(\varphi^\alpha_u(y)) \psi_u$, with $\psi_0 = 1$. By uniqueness, \eqref{eq:derivative-flow} follows.

Therefore, the function $y \mapsto H^y_\alpha(t)$ is $C^1$ with
\[ \frac{\d}{\d y}H^y_\alpha(t) = -H^y_\alpha(t)\int_0^t{f'(\varphi^\alpha_u(y)) \frac{\d}{\d y} \varphi^\alpha_u(y) \d u} = - H^y_\alpha (t)  \frac{f(\varphi^\alpha_t(y)) - f(y)}{b(y) + \alpha}.  \]
We used \eqref{eq:derivative-flow} to obtain the last equality.
Let $A(t)$ be equal to the right-end side of \eqref{eq:formula-for-Psi}.
Plugging the explicit expression of $\nu^\infty_\alpha$ (see~\eqref{eq:invariant-measure}) and using that $J = \alpha/\gamma(\alpha)$, we find
\begin{align*}
	A(t) &= \alpha \int_0^{\sigma_\alpha}{ H^x_\alpha(t)\frac{f(\varphi^\alpha_t(x)) - f(x)}{(b(x) + \alpha)^2} \exp\left(-\int_0^x{\frac{f(y)}{b(y) + \alpha} \d y} \right) \d x } \\
	     &=  \alpha \int_0^\infty{ \exp \left(-\int_0^t{ f(\varphi^\alpha_{s + u}(0)) \d s} \right)  \frac{f(\varphi^\alpha_{t+u}(0)) - f(\varphi^\alpha_u(0))}{b(\varphi^\alpha_u(0)) + \alpha}  H_\alpha(u)\d u }.
\end{align*}
To obtain the last equality we made the change of variables $x = \varphi^{\alpha}_u(0)$ and then $y = \varphi^\alpha_\theta(0)$.
Hence, we find that $A(t) = \Psi_\alpha(t)$ as claimed.  
\end{proof}
Then, we define
\[ \Theta_\alpha(t) := J \int_{\R_+} \frac{\d}{\d y} \E_y f(Y^\alpha_{t}) \nu^\alpha_\infty(\d y). \]
Recall that $\lambda_\alpha > 0$ is defined in Assumption~\ref{ass:final} and that $D_\alpha := \{z \in \C; ~\Re(z) > -\lambda_\alpha \}$.  
\begin{lemma}
	For all $z \in D_\alpha$, it holds that $\hat{\Psi}_\alpha(z) = \hat{H}_\alpha(z)$ if and only if $\hat{\Theta}_\alpha(z) = 1$. 
\end{lemma}
\begin{proof}
	Recall that $r^x_\alpha(t) := \E_x f(Y^\alpha_t)$. Using Lemma~\ref{lem:volterra-eq}, we have
	\[ r^x_\alpha = K^x_\alpha + r_\alpha * K^x_\alpha. \]
	We used here the notation $(r_\alpha * K^x_\alpha)(t) := \int_0^t{ r_\alpha(t-s) K^x_\alpha(s) \d s}$. 
	We differentiate with respect to $x$ and obtain
	\[ J \frac{\d }{\d x} r^x_\alpha = J \frac{\d }{\d x} K^x_\alpha + r_\alpha * \left[ J \frac{\d }{\d x} K^x_\alpha \right] . \]
Let
	\begin{equation}
		\label{eq:def-of-Xi}
	\Xi_\alpha(t) := \frac{\d }{\d t} \Psi_\alpha(t).
\end{equation}
In view of \eqref{eq:formula-for-Psi} and $\frac{\d }{\d t} H^{x}_\alpha(t) = - K^x_\alpha(t)$, we have:
\[ 	\Xi_\alpha(t) =  J \int_{\R_+} \frac{\d }{\d y} K^y_\alpha(t) \nu^\alpha_\infty(\d y). \]
We deduce that
	\[ \Theta_\alpha = \Xi_\alpha + r_\alpha * \Xi_\alpha. \]
	Taking the Laplace transform, we find that for all $z \in \C$ with  $\Re(z) > 0$:
\[ \widehat{\Theta}_\alpha(z) = \widehat{\Xi}_\alpha(z) + \widehat{r}_\alpha(z) \widehat{\Xi}_\alpha(z). \]
	Therefore, we have
	\begin{align}
	\widehat{\Theta}_\alpha(z) &= \widehat{\Xi}_\alpha(z) \left[ 1 + \widehat{r}_\alpha (z) \right] \nonumber \\
				   &= \widehat{\Xi}_\alpha(z) \left[ 1 + \frac{\widehat{K}_\alpha(z) }{1 - \widehat{K}_\alpha (z)} \right] \quad (\text{using } r_\alpha = K_\alpha + K_\alpha * r_\alpha)  \nonumber \\
&=  \frac{ \widehat{\Xi}_\alpha(z)   }{z \widehat{H}_\alpha (z)} \quad (\text{using } \widehat{K}_\alpha(z) = 1 - z \widehat{H}_\alpha(z)) \label{eq:alternative formula for Laplace Theta} \\
&= \frac{ \widehat{\Psi}_\alpha (z) }{ \widehat{H}_\alpha(z) }  \quad (\text{using } \Psi_\alpha(0) = 0 \text{ and } \frac{\d}{\d t} \Psi_\alpha = \Xi_\alpha). \nonumber 
\end{align}
Because the left-hand side and the right-hand side are two analytic functions on $D_\alpha$, the 
equality is in fact valid on $D_\alpha$ and so 
\[ \forall z \in D_\alpha, \quad \widehat{\Theta}_\alpha(z) = 1 \iff \widehat{\Psi}_\alpha(z) = \widehat{H}_\alpha(z). \]
\end{proof}
We finally consider $\Omega_\alpha(t)$ the solution of the Volterra integral equation
\[ \forall t \geq 0, \quad \Omega_\alpha(t) = \Theta_\alpha(t) + \int_0^t{ \Omega_\alpha(t-s) \Theta_\alpha(s) \d s}. \]
\begin{remark}
	The function $\Omega_\alpha(t)$ has a simple probabilistic interpretation using the McKean-Vlasov equation~\eqref{eq:McKeanVlasov}. For all $\epsilon > 0$, let 
	$\nu_\epsilon := \Law(X^{\nu_\infty}_0 + \epsilon)$, with $\Law(X^{\nu_\infty}_0) = \nu_\infty$. Then
	\[ \Omega_\alpha(t) = \lim_{\epsilon \downarrow 0} \frac{\E f(X^{\nu_\epsilon}_t) - \E f(X^{\nu_\infty}_t)}{\epsilon}. \]
	Similarly, it holds that
	\[ \Theta_\alpha(t) = \lim_{\epsilon \downarrow 0} \frac{\E f(Y^{\alpha, \nu_\epsilon}_{t}) - \E f(Y^{\alpha, \nu_\infty}_{t})}{\epsilon}, \]
	where $(Y^{\alpha, \nu}_t)$ is the solution of \eqref{eq:linear-process}.
	We refer to \cite{cormier2023stability} for these probabilistic interpretations as well as the connection with Lions derivatives.
\end{remark}
\begin{lemma}
	\label{lem:Omega-t-lambda}
	For all $\lambda < \lambda'_\alpha$, where $\lambda'_\alpha$ is given by \eqref{eq:spectral-assumption}, we have $\sup_{t \geq 0} \abs{\Omega_\alpha(t)} e^{\lambda t} < \infty$.
\end{lemma}
In other words, $\lambda'_\alpha$ gives the rate of convergence of $\Omega_\alpha(t)$ towards zero.
\begin{proof}
	Let $\lambda < \lambda'$, $K_t := e^{\lambda t} \Theta_\alpha(t)$ and $R_t := e^{\lambda t} \Omega_\alpha(t)$. It holds that $K \in L^1(\R_+)$.
	By assumption, it holds that $\hat{K}(z) \neq 0$ for all $\Re(z) \geq 0$. 
	Therefore, \cite[Ch. 2, Th. 4.1]{gripenberg_volterra_1990} applies and so $R \in L^1(\R_+)$. 
	Finally, using that $R = K + K*R$ we find that $R$ is also bounded.
\end{proof}
\subsection{Proof of Theorem~\ref{th:stability}}
\subsubsection*{A sensitivity formula}
Following~\cite{cormier2023stability}, we first show the following ``sensitivity'' formula:
\begin{proposition} \label{prop:sensitivity}
	Let $k \in C([0, t]; \R)$ and $\alpha \in \R_+$, such that $\inf_{s \in [0, t]} (\alpha + k_s) \geq 0$. Provided that $\Law(Y^{\alpha+k}_0) = \Law(Y^\alpha_0)$, it holds that
	\[  \E g(Y^{\alpha+k}_t) - \E g(Y^\alpha_t) = \int_0^t \int_{\R_+} \left[ \frac{\d}{\d y} \E_y g(Y^\alpha_{t-\theta}) \right]  k_\theta \Law(Y^{\alpha+k}_{\theta})(\d y) \d \theta. \]
\end{proposition}
\begin{proof}
	This is a Trotter-Kato formula. Define for all $s \in (0, t]$ and for all $y \in \R_+$:
	\[ \phi(s, y) := \E_y g(Y^\alpha_{t-s}). \]
	The function $\phi$ is $C^1_b(\R_+ \times \R_+)$ and
	\[ \frac{\partial}{\partial s} \phi(s, y) =  -\cL^\alpha \phi(s, y), \]
	where $\cL^\alpha$ is the generator of $(Y^\alpha_t)$, solution of \eqref{eq:linear-process}. This generator acts on the marginal function $\phi(s, \cdot)$ and is given by
	\[ \forall g \in C^1_b(\R_+), \quad \cL^\alpha g (y) := g'(y) (b(y) + \alpha) + (g(0) - g(y)) f(y). \]
	In addition, the generator of $(Y^{\alpha+k}_s)$ satisfies $(\cL^{\alpha+k}_s - \cL^\alpha) g(y) = g'(y) k_s$.
	Therefore, by Ito's formula, we obtain:
	\[ \E \phi(s, Y^{\alpha+k}_s) = \E \phi(0, Y^{\alpha+k}_0) + \E \int_0^s{ \frac{\partial}{\partial y} \phi(u, Y^{\alpha+k}_u) k_u \d u  }. \]
	Replacing $\phi$ by its definition, we deduce that
	\[ \E \phi(s, Y^{\alpha+k}_s) = \E \phi(0, Y^{\alpha+k}_0) + \E \int_0^s{ \int_{\R_+} \left[ \frac{\d}{\d y} \E_y g(Y^\alpha_{t-u})  \right] k_u \Law(Y^{\alpha+k}_u)(\d y) \d u  }. \] 
	Finally, we let $s$ converge to $t$. Using the Markov property at time $s = 0$ and the fact that $\phi(t, y) = g(y)$, we find that the stated formula holds.
\end{proof}
\begin{corollary} \label{cor:control-distance}
	It holds that
	\[ \norm{\Law(Y^{\alpha+k, \nu}_t) - \Law(Y^{\alpha, \nu}_t)}_0 \leq C_* \int_0^t{ e^{-\lambda_* (t-s) } \abs{k_s} \d s}.   \]
\end{corollary}
\begin{proof}
	Let $g \in \Lip_1(\R_+)$ with $\norm{g}_\infty \leq 1$. We have by Proposition~\ref{prop:sensitivity}:
\[  \abs{\E g(Y^{\alpha+k, \nu}_t) - \E g(Y^{\alpha, \nu}_t)} \leq \int_0^t \left\{ \sup_{y \in \R_+} \left| \frac{\d}{\d y} \E_y g(Y^\alpha_{t-s}) \right| \right\} \abs{k_s} \d s. \]
By Proposition~\ref{prop:estimate-linear-process}, we have for $y \neq y'$:
\[ \abs{ \E_y g(Y^\alpha_t) - \E_{y'} g(Y^\alpha_t) } \leq C_* e^{-\lambda_* t} \norm{\delta_{y} - \delta_{y'}}_0.  \]
As  $\norm{\delta_{y} - \delta_{y'}}_0 = |y-y'|$ for $|y-y'| \leq 1$, we deduce that $\left| \frac{\d}{\d y} \E_y g(Y^\alpha_{t-s}) \right| \leq C_* e^{-\lambda_*(t-s)}$. \end{proof}
\subsubsection*{Control of the non-linear interactions}
We define:
\begin{align*}
	\varphi^\nu_t &:= J \E f(Y^{\alpha, \nu}_t) - \alpha \\
	k^\nu_t &= J \E f(X^\nu_t) - \alpha.
\end{align*}
We prove that:
\begin{proposition}
For all $T > 0$, there is a constant $C_T$ such that for all $t \in [0, T]$ and for all $\nu \in \cP(\R^d)$:
\begin{enumerate}
	\item $\abs{k^\nu_t} \leq C_T \norm{\nu - \nu_\infty}_0$. 
	\item $\left|k^\nu_t - \varphi^\nu_t - \int_0^t{\Theta_\alpha(t-s) k^\nu_s \d s} \right| \leq  C_T \left(\norm{\nu - \nu_\infty}_0 \right)^2$.
	\item $\left|k^\nu_t - \varphi^\nu_t - \int_0^t{\Omega_\alpha(t-s) \varphi^\nu_s \d s} \right| \leq  C_T  \left( \norm{\nu - \nu_\infty}_0 \right)^2$. 
\end{enumerate}
\end{proposition}
\begin{proof}
	The first point is a consequence of $\norm{f}_\infty + \norm{f'}_\infty < \infty$ and of Theorem~\ref{th:stab}.
	For the second point, we note that $\E f(X^\nu_t) = \E f(Y^{\alpha + k^\nu, \nu}_t)$.
	We define $\psi_t(y) = \frac{\d}{\d y} \E_y f(Y^\alpha_t)$. We have, using Proposition~\ref{prop:sensitivity} with $g = J f$:
	\begin{align*}
		k^\nu_t - \varphi^\nu_t &= J \E f(Y^{\alpha + k^\nu, \nu}_t) - J \E f(Y^{\alpha, \nu}_t) \\ 
					&= J\int_0^t \int_{\R_+} \frac{\d }{\d y} \E_y f(Y^\alpha_{t-s}) k^\nu_s  \Law(Y^{\alpha+k^\nu, \nu}_s)(\d y) \d s \\
				     &= J\int_0^t \E \psi_{t-s}(Y^{\alpha + k^\nu, \nu}_s)  k^\nu_s \d s \\
				     &= \int_0^t \Theta_\alpha(t-s) k^\nu_s \d s  + J\int_0^t \left[ \E \psi_{t-s}(Y^{\alpha + k^\nu, \nu}_s)  - \E \psi_{t-s}(Y^{\alpha, \nu_\infty}_s) \right] k^\nu_s \d s. 
	\end{align*}
	Because $f$ and $b$ are assumed to be $C^2$, there exists a constant $C_T$ such that
	\[ \forall t \in [0, T], \quad  \norm{\psi_t}_\infty + \norm{\partial_y \psi_t}_\infty \leq C_T. \]
	Therefore, by Theorem~\ref{th:stab} we have:
	\begin{align*}
	\left| \E \psi_{t-s}(Y^{\alpha + k^\nu, \nu}_s)  - \E \psi_{t-s}(Y^{\alpha, \nu_\infty}_s) \right| &=  \left| \E \psi_{t-s}(X^{\nu}_s)  - \E \psi_{t-s}(X^{\nu_\infty}_s) \right|  \leq C_T \norm{\nu - \nu_\infty}_0.
	\end{align*}
	Using the first point, we obtain the stated inequality. The last point is obtained by iterating the estimate of the second point, as in \cite{cormier2023stability}.
\end{proof}
Exactly as in \cite[Lem. 2.20]{cormier2023stability}, we deduce from this last result and from Corollary~\ref{cor:control-distance} that:
\begin{lemma} \label{lem:key-lemma}
	Let $\lambda \in (0, \lambda'_\alpha)$. There exists a constant $C_\lambda$ such that for all $T > 0$, there exists $C_T > 0$: for all $\nu \in \cP(\R)$, for all $t \in [0, T]$,
	\[ \norm{\Law(X^\nu_t) -  \nu_\infty}_0 \leq C_\lambda e^{ - \lambda t} \norm{\nu-\nu_\infty}_0 + C_T \left( \norm{\nu-\nu_\infty}_0 \right)^2. \]
\end{lemma}
We used crucially here that $\sup_{t \geq 0} \abs{\Omega_\alpha(t)} e^{\lambda t} < \infty$, see Lemma~\ref{lem:Omega-t-lambda}.
The proof of Theorem~\ref{th:stability} is easily deduced from Lemma~\ref{lem:key-lemma}, exactly as in \cite{cormier2023stability}.
\subsection{Proof of Proposition~\ref{prop:structural-stability}}
The first point is to note that under the assumption $\inf_{[0, \sigma_\alpha)} f + b' \geq 0$, we can integrate by parts $\Psi_\alpha$ and 
$\Xi_\alpha$:
\begin{lemma} \label{lem:alternative-formula-for-Psi-and-Xi} We have
	\begin{enumerate}
		\item The following limit exists and is finite: $\nu^\infty_\alpha(\sigma_\alpha) := \lim_{x \uparrow \sigma_\alpha}{\nu^\infty_\alpha(x)}  < \infty$. 
		\item Define $C_\alpha := \frac{b(0) + \alpha}{\gamma(\alpha)} \nu^\infty_\alpha(\sigma_\alpha) $ 
		and 
		\begin{equation}
		\Upsilon_\alpha(t) := C_\alpha H^{\sigma_\alpha}_\alpha(t) +  \int_0^\infty{ H_\alpha(t+u) \left[ 
		f(\varphi^\alpha_u(0)) + b'(\varphi^\alpha_u(0)) \right] \frac{b(0) + \alpha}{ b(\varphi^\alpha_u(0) ) 
		+ \alpha}  \d u}. 
		\label{eq:formule donnant Upsilon lorsque f + b' assymptotiquement strictement positif}
		\end{equation}
		It holds that for all $t \geq 0$:
		\begin{equation}
			\label{eq:alternative formule for Psi alpha}
			 \Psi_\alpha(t) = \frac{\alpha}{b(0) + \alpha} \left[  H_\alpha (t) - 
		\Upsilon_\alpha (t) 
		\right]. 
		\end{equation}
		\item Define $\Lambda_\alpha(t) := -\frac{\d }{\d t} \Upsilon_\alpha(t)$. One has for all $t \geq 0$
		\begin{equation}
		\Lambda_\alpha(t) = C_\alpha K^{\sigma_\alpha}_\alpha(t) + \int_0^\infty{ K_\alpha(t+u) \left[ 
		f(\varphi^\alpha_u(0)) + b'(\varphi^\alpha_u(0)) \right] \frac{b(0) + \alpha}{ b(\varphi^\alpha_u(0) ) 
		+ \alpha}  \d u}. 
		\label{eq:formule donnant Lambda lorsque f + b' assymptotiquement strictement positif}
		\end{equation}
		Moreover, for all $t \geq 0$
		\begin{equation}\Xi_\alpha(t) = \frac{\alpha}{b(0) + \alpha} \left[ \Lambda_\alpha (t) - 
		K_\alpha (t) \right]. 
		\label{eq:alternative formule for Xi alpha}
		\end{equation}
	\end{enumerate}
\end{lemma}
\begin{proof}[Proof of Lemma~\ref{lem:alternative-formula-for-Psi-and-Xi}]
	To prove the first point, we use the explicit formula of the invariant measure \eqref{eq:invariant-measure}: we find that for all $x < \sigma_\alpha$
	\[ \frac{\d}{\d x} \nu^\alpha_\infty(x) = -\gamma(\alpha) \frac{f(x) + b'(x)}{(b(x) + \alpha)^2} \exp \left(-\int_0^x \frac{f(y)}{b(y) + \alpha} \d y \right) \leq 0. \] 
	Therefore, $x \mapsto \nu^\alpha_\infty(x)$ is non-increasing and so $ \lim_{x \uparrow \sigma_\alpha}{\nu_\infty^\alpha(x)}$ exists and is finite (it is equals to zero if $\sigma_\alpha = \infty$, and might be non-null in the case where $\sigma_\alpha < \infty$).  
	To prove the second point, we integrate by parts the right-hand side of \eqref{eq:formula-for-Psi} and find
	\[ \Psi_\alpha(t) = \frac{\alpha}{b(0) + \alpha} \left[ H_\alpha(t) - C_\alpha 
	H^{\sigma_\alpha}_\alpha(t)\right] + J \int_0^{\sigma_\alpha}{ H^x_\alpha(t) \frac{\d}{\d x} 
	\nu_\infty^\alpha(x) \d x }. \]
	The last term is equal to:
	\[
	J \int_0^{\sigma_\alpha}{ H^x_\alpha(t) \frac{\d}{\d x} \nu_\infty^\alpha(x) \d x } = - \alpha \int_0^{\sigma_\alpha} H^x_\alpha(t)  \frac{f(x) + b'(x)}{(b(x) + \alpha)^2} \exp \left(-\int_0^x \frac{f(y)}{b(y) + \alpha} \d y \right) \d y. \] 
	We make the changes of variables $y = \varphi^\alpha_\theta(0)$ and $x = \varphi^\alpha_u(0)$ 
	and obtain
	\[ \Psi_\alpha(t) = \frac{\alpha}{b(0) + \alpha} \left[ H_\alpha(t) - C_\alpha 
	H^{\sigma_\alpha}_\alpha(t) \right]  - \alpha \int_0^{\infty}{ 
	H^{\varphi^\alpha_u(0)}_\alpha(t) \frac{f(\varphi^\alpha_u(0))  + 
	b'(\varphi^\alpha_u(0))}{b(\varphi^\alpha_u(0)) + \alpha} H_\alpha(u) \d u }.   \] 
	Using that $H^{\varphi^\alpha_u(0)}_\alpha(t) H_\alpha(u) = H_\alpha(t+u)$, we obtain the stated 
	formula. Finally, recall that $\Xi_\alpha(t) = \frac{\d}{\d t} \Psi_\alpha(t)$. Therefore, the third point is obtained by differentiating the second point with respect to $t$. 
\end{proof}
\begin{proof}[Proof of Proposition~\ref{prop:structural-stability}, first point]
	Recall that $\Xi_\alpha(t)$ is given by \eqref{eq:def-of-Xi} and satisfies for all $z \in D_\alpha$, $\widehat{\Xi}_\alpha(z) = z \widehat{\Psi}_\alpha(z)$. 
	Similarly, it holds that $\widehat{K}_\alpha(z) = 1 - z \widehat{H}_\alpha(z)$.	
	Let $z_* \in D_\alpha$ such that $\widehat{H}_\alpha(z_*) = \widehat{\Psi}_\alpha(z_*)$. We deduce that
	\[  1 - \widehat{K}_\alpha(z_*) = \widehat{\Xi}_\alpha(z_*). \]
	Using Lemma~\ref{lem:alternative-formula-for-Psi-and-Xi}, we have
	\begin{equation} \label{eq:reduced equation CS for constant b} b(0) + \alpha = b(0) \widehat{K}_\alpha(z_*) + \alpha \widehat{\Lambda}_\alpha(z_*). \end{equation}
	We show that $\Re(z_*) < 0$. Indeed, we have
	\[ \Re(z) > 0 \implies  |b(0) \widehat{K}_\alpha (z) + \alpha \widehat{\Lambda}_\alpha(z) | < b(0) 
	|\widehat{K}_\alpha (z) | + \alpha |\widehat{\Lambda}_\alpha(z)| < b(0) + \alpha,\]
	and so necessarily, $\Re(z_*) \leq 0$. 
	We used here that both $K_\alpha(t)$ and $\Lambda_\alpha(t)$ are the densities of probability measures.
	In addition, if $z_* = i \omega$ for some $\omega > 0$, then 

	\[ \Re \left[ b(0)(1 - \widehat{K}_\alpha(i \omega) ) + \alpha( 1 - \widehat{\Lambda}_\alpha (i \omega) ) \right] = 
	\int_0^\infty{ [1 - \cos(\omega t)] (b(0) K_\alpha(t) + \alpha \Lambda_\alpha(t) ) \d t}.  \]
	Because for $t\in\mathbb{R}_+$, $1 - \cos(\omega t) > 0$ almost everywhere, the right-hand side is null 
	only if almost everywhere
	\[ b(0) K_\alpha(t) + \alpha \Lambda_\alpha(t) = 0. \] 
	This leads to a contradiction because by \eqref{eq:formule donnant Lambda lorsque f + b' assymptotiquement strictement positif}, we have $\Lambda_\alpha(t) \geq 0$. In addition, $K_\alpha(t) \geq 0$ and the total mass of $K_\alpha$ is one.
	Altogether, we have proved that
	\[ \forall z_* \in D_\alpha, \quad \widehat{H}_\alpha(z_*) = \widehat{\Psi}_\alpha(z_*) \implies \Re(z_*) < 0. \]
	Using the Riemann–Lebesgue lemma, we deduce that $\lambda'_\alpha > 0$.
\end{proof}
\begin{proof}[Proof of Proposition~\ref{prop:structural-stability}, second point]
	Assume that $\inf_{x \geq 0} f(x) + b'(x) \geq 0 $ and $b(0) >0$.
	The number of invariant probability measures of \eqref{eq:McKeanVlasov} is given by the number of solutions of the equation $\alpha = J \gamma(\alpha), \alpha \geq 0$.
	We first prove that the continuous function 
	$G: \alpha \mapsto 
	\frac{\alpha}{\gamma(\alpha)}$
	is strictly increasing on $\mathbb{R}_+$. Note first the identity:
	\[ \forall t \geq 0, \quad \left[ b(\varphi^\alpha_t(0)) + \alpha \right] \exp \left(- 
		\int_0^t{b'(\varphi^\alpha_u(0)) \d u } \right) = b(0) + \alpha.\]
	We deduce that for all $\alpha > 0$
	\begin{align*}
	G(\alpha) = \frac{\alpha}{\gamma(\alpha)} & = \alpha 
	\int_0^\infty{ H_\alpha(t) \d t} \\
	&=\frac{\alpha}{b(0) + \alpha } \int_0^\infty{  \left[ b(\varphi^\alpha_t(0)) + \alpha \right] \exp \left(- 
	\int_0^t{b'(\varphi^\alpha_u(0)) \d u } \right) H_\alpha(t) \d t  } \\
	&= \frac{\alpha}{b(0) + \alpha } \int_0^\infty{  \left[ b(\varphi^\alpha_t(0)) + \alpha \right] 
	\exp \left(-\int_0^t{ (f + b') (\varphi^\alpha_u(0)) \d u} \right) \d t  }.
	\end{align*}
	The changes of variables $y = \varphi^\alpha_u(0)$ and $x = \varphi^\alpha_t(0)$ show that
	\[ \frac{\alpha}{\gamma(\alpha)}  = \frac{\alpha}{b(0) + \alpha } \int_0^{\sigma_\alpha}{ \exp\left(- 
	\int_0^x{ \frac{(f + b')(y)}{b(y) + \alpha } \d y} \right)  \d x}. \]
	Note that the function $\alpha \mapsto \frac{\alpha}{b(0) + \alpha }$ is non-decreasing and $\alpha 
	\mapsto \sigma_\alpha$  is strictly increasing.
	Moreover, because $f + b' \geq 0$, for all fixed $x$, the function
	\[ \alpha \mapsto  \exp \left(- \int_0^x{ \frac{(f + b')(y)}{b(y) + \alpha } \d y} \right) \]
	is non-decreasing. 
	So $G$ is strictly increasing. Because $b(0) > 0$, we have $\gamma(0) > 0$ and so $G(0) = 0$. Therefore, for all $J \geq 0$, the equation $G(\alpha) = J$ has a unique solution. \end{proof}
\section{An illustrated example}
\label{sec:example}
To illustrate the results, we consider
\[ f(x) = x^2 \quad \text{ and } \quad  b(x) = - x, \quad \forall x \geq 0. \]
There is a slight technical difficulty: $f$ and $f'$ are not bounded and so we cannot directly apply our results. For $A > 0$, we denote by
\[ \cM_A := \{ \nu \in \cP(\R_+); \quad J \langle f, \nu \rangle \leq A, \quad  \nu([0, A]) = 1 \}. \]
\begin{lemma}
	\label{lem:suppA}
	Let $J \geq 0$. There exists a constant $A > 0$ large enough such that for any initial condition $\nu \in \cM_A$, there is a unique path-wise solution to \eqref{eq:McKeanVlasov} and for all $t \geq 0$, $\Law(X^\nu_t) \in \cM_A$. 
\end{lemma}
\begin{proof}
	Existence and uniqueness of the solution of \eqref{eq:McKeanVlasov} is not problematic since the initial condition is compactly supported, see \cite{fournier_toy_2016}.
	The existence of the constant $A$ is shown in \cite{CTV}. The idea is that by Ito's formula, 
	\begin{align*} \frac{\d }{\d t}  \E f(X^\nu_t) & = \E f'(X^\nu_t) [b(X^\nu_t) + J \E f(X^\nu_t)] - \E  f^2(X^\nu_t) \\
		& \leq \frac{C^2}{2} - \frac{1}{2}  \E f^2(X^\nu_t)\\
		& \leq \frac{C^2}{2} - \frac{1}{2} (\E f(X^\nu_t))^2, 
	\end{align*}
for some constant $C$ only depending on $J$. We used the Cauchy-Schwarz inequality to obtain the last estimate. Then we deduce that 
$\E f(X^\nu_0) \leq C$ implies $\E f(X^\nu_t) \leq C$ for all $t \geq 0$. Therefore, we can choose $A = J C$ to obtain the result.
	\end{proof}
	We equip $\cM_A$ with the distance $\norm{\nu - \mu}_0$ (we could also use here the standard Wasserstein distance $W_1(\nu, \mu)$: $W_1$ and $\norm{\cdot}_0$ are equivalent on $\cM_A$, as the probability measures are compactly supported). 
	In view of Lemma~\ref{lem:suppA}, as soon as the initial condition belongs to $\cM_A$, everything happens as if $f$ and $f'$ were bounded on $\R_+$. Therefore Theorem~\ref{th:stab} and \ref{th:stability} holds with $\cP(\R_+)$ being replaced with $\cM_A$.

	Our goal is now to describe all the invariant probability measures and to find their stability properties. The first step to apply Theorem~\ref{th:stability} is to specify the values of $\widehat{H}_\alpha$ and $\widehat{\Psi}_\alpha$. 
	\subsection*{Computation of $\widehat{H_\alpha}(z)$ and $\widehat{\Psi_\alpha}(z)$}
	Let $\alpha > 0$. Recall that $H_\alpha(t)$ is defined by \eqref{eq:def-of-Halpha} and that $\Psi_\alpha(t)$ is defined by \eqref{eq:def-of-Psi-alpha}. Let: 
	\[ \forall \theta \in [0, 1), \quad w(\theta) := \theta + \frac{\theta^2}{2} + \log(1-\theta) = - \sum_{k \geq 3} \frac{\theta^k}{k}. \]
	In this section, we prove that:
	\begin{proposition}
		\label{prop:expression-of-Psi-and-H-hat}
		Let $\psi(z, x) := \frac{2-2(1-x)^z - 2xz - (1-x)^z x^2 (1-z) z }{(1-x)(1-z)z(1+z)}$.
		It holds that for all $\Re(z) > -\alpha^2$:
		\[ \widehat{H_\alpha}(z) = \int_0^1 (1-x)^{z-1} e^{\alpha^2 w(x)} \d x \]
		and
		\[ \widehat{\Psi_\alpha}(z) = \alpha^2 \int_0^1{ \psi(z, x) e^{\alpha^2 w(x)} \d x}. \]
	\end{proposition}
	To proceed, we introduce some notations. First, recall that $K_\alpha(t) = - \frac{\d}{\d t} H_\alpha(t)$.
	We also consider
	\begin{align*} \forall \theta \in [0, 1), \quad \widehat{H^{[\theta]}_\alpha}(z) &:= \int_0^\infty e^{-zt} H_\alpha(t - \log(1-\theta)) \d t, \\
 \widehat{K^{[\theta]}_\alpha}(z) &:= \int_0^\infty e^{-zt} K_\alpha(t - \log(1-\theta)) \d t. 
	\end{align*}
	\begin{lemma}
		\label{lem:formula for H-K-theta-z}
	The following identities hold:
		\begin{align*}
			\widehat{H^{[\theta]}_\alpha}(z) &=  (1-\theta)^{-z} \int_\theta^1 (1-x)^{z-1} e^{\alpha^2 w(x) } \d x,  \\
			\widehat{K^{[\theta]}_\alpha}(z) &= e^{\alpha^2 w(\theta)} - z \widehat{H^{[\theta]}_\alpha}(z). 
		\end{align*}
		In addition, we have
		\[  e^{\alpha^2 w(\theta)} = (z+\alpha^2) \widehat{H^{[\theta]}_\alpha}(z) - \alpha^2 \int_\theta^1 \left(\frac{1-x}{1-\theta}\right)^z (1+x) e^{\alpha^2 w(x)} \d x.  \]
	\end{lemma}
	\begin{proof}
		We have $\varphi^\alpha_t(0) = \alpha(1-e^{-t})$ and so
		\[ H_\alpha(t) = \exp \left( -\int_0^t{ \alpha^2 (1-e^{-u}})^2 \d u \right) = \exp \left(\alpha^2[\frac{3}{2} + \frac{e^{-2 t}}{2} - 2 e^{-t} - t] \right). \]
		Therefore, we find that
		\[ \forall x \in [0, 1), \quad H_\alpha(-\log(1-x)) = e^{\alpha^2 w(x)}. \]
		So, the change of variable $x = 1-e^{-t}$ shows that
		\begin{align*}
			\widehat{H^{[\theta]}_\alpha}(z) & = (1-\theta)^{-z} \int_{-\log(1-\theta)}^\infty e^{-zt} H_\alpha(t) \d t \\
							 &= (1-\theta)^{-z }\int_\theta^1 (1-x)^{z-1} e^{\alpha^2 w(x) } \d x. 
		\end{align*}
		This proves the first equality.	The second equality is obtained by an integration by parts, using that $\frac{\d}{\d t} H_\alpha(t - \log(1-\theta)) = - K_\alpha(t - \log(1-\theta))$. To obtain the last identity, we note that:
		\[  - \frac{\d }{\d x} \left[ (1-x)^z e^{\alpha^2 w(x)} \right] = (z+\alpha^2) (1-x)^{z-1} e^{\alpha^2 w(x)} - \alpha^2 (1-x)^z(1+x) e^{\alpha^2 w(x)}.   \]
		We then integrate this equality from $x = \theta$ to $x = 1$. This gives the result.
	\end{proof}
	Using this Lemma, we can finally deduce the expression of $\widehat{\Psi_\alpha}$:
	\begin{proof}[Proof of Proposition~\ref{prop:expression-of-Psi-and-H-hat}]
		For $u \in \R_+$, let $\theta(u) := 1-e^{-u}$. Recall that $\Psi_\alpha(t)$ is given by \eqref{eq:def-of-Psi-alpha}. Therefore, it holds that
		\begin{align*}
			\widehat{\Psi_\alpha}(z) &= \alpha \int_0^\infty \frac{1}{b(\varphi^\alpha_u(0)) + \alpha} \left[ \int_0^\infty e^{-zt} K_\alpha(t+u) \d t - f(\varphi^\alpha_u(0))  \int_0^\infty e^{-zt} H_\alpha(t+u) \d t \right] \d u \\
						 &= \int_0^\infty{ e^{u} \left[ \widehat{K^{[\theta(u)]}_\alpha}(z) - \alpha^2 (\theta(u))^2 \widehat{H^{[\theta(u)]}_\alpha}(z) \right] \d u }  \\
						 &= \int_0^1 (1-\theta)^{-2} \left[ \widehat{K^{[\theta]}_\alpha}(z) - \alpha^2 \theta^2 \widehat{H^{[\theta]}_\alpha}(z) \right] \d \theta.  
		\end{align*}
	We made the change of variable $\theta = 1-e^{-u}$.
	We now use Lemma~\ref{lem:formula for H-K-theta-z} and find that
	\begin{align*} \widehat{K^{[\theta]}_\alpha}(z) - \alpha^2 \theta^2 \widehat{H^{[\theta]}_\alpha}(z) &= \alpha^2(1-\theta^2) \widehat{H^{[\theta]}_\alpha}(z) - \alpha^2  \int_\theta^1 \left(\frac{1-x}{1-\theta}\right)^z (1+x) e^{\alpha^2 w(x)} \d x.   \\
		&=\int_\theta^1 \left[  \alpha^2(1-\theta^2)  \frac{(1-x)^{z-1}}{(1-\theta)^z} - \alpha^2 \left(\frac{1-x}{1-\theta}\right)^z (1+x)  \right] e^{\alpha^2 w(x) } \d x. 
	\end{align*}
	Therefore, we find that 
		\begin{align*} \widehat{\Psi_\alpha}(z) &= \alpha^2 \int_0^1 \frac{1 - \theta^2 }{(1-\theta)^z} \int_\theta^1 (1-x)^{z-1} e^{\alpha^2 w(x)} \d x \d \theta \\
	&\quad - \alpha^2 \int_0^1 \frac{1}{(1-\theta)^z} \int_\theta^1 (1-x)^z (1+x) e^{\alpha^2 w(x)} \d x \d \theta.  \end{align*}
	To obtain the stated result, it suffices to integrate by parts this equality.
	\end{proof}
\subsection{Description of the invariant probability measures}
The following proposition gives the number of invariant measures of the non-linear 
equation~\eqref{eq:McKeanVlasov}. This result is conjectured to be true in \cite[Section 
7.2.3]{zbMATH06706511}.
\begin{figure}[ht]
		\centering
\begin{tikzpicture}[scale=0.8]
\begin{axis}[xlabel={$\alpha$}, ylabel={$J(\alpha) = \frac{\alpha}{\gamma(\alpha)}$}, 
title={Graph of the function $\alpha \mapsto J(\alpha) := \frac{\alpha}{\gamma(\alpha)}$. }, 
	 grid=major, grid style=dashed, xmajorgrids, ymajorgrids]
	\addplot[no marks, blue, very thick]
	coordinates {
		(0.3,3.7676033626182988)
		(0.39591836734693875,3.0852009670878378)
		(0.49183673469387756,2.7086065212296466)
		(0.5877551020408164,2.483311659752082)
		(0.6836734693877551,2.3419060086028693)
		(0.7795918367346939,2.2506322736080993)
		(0.8755102040816326,2.190970380161363)
		(0.9714285714285714,2.1520789270667526)
		(1.0673469387755101,2.1272923668818198)
		(1.163265306122449,2.112348901091792)
		(1.2591836734693878,2.1044333705836413)
		(1.3551020408163266,2.1016335501065315)
		(1.4510204081632654,2.102618496022046)
		(1.546938775510204,2.106441710972514)
		(1.6428571428571428,2.11241712528422)
		(1.7387755102040816,2.120038919133254)
		(1.8346938775510204,2.1289284757068847)
		(1.9306122448979592,2.138798540556206)
		(2.026530612244898,2.1494285380396523)
		(2.122448979591837,2.1606472704003914)
		(2.2183673469387757,2.172320594299115)
		(2.3142857142857145,2.1843425114303074)
		(2.410204081632653,2.1966286386590106)
		(2.5061224489795917,2.2091113613854234)
		(2.6020408163265305,2.221736194092055)
		(2.6979591836734693,2.234459017807991)
		(2.793877551020408,2.247243962207478)
		(2.889795918367347,2.260061766870262)
		(2.9857142857142858,2.272888502405327)
		(3.0816326530612246,2.285704564459513)
		(3.1775510204081634,2.298493876526521)
		(3.273469387755102,2.311243253872034)
		(3.369387755102041,2.3239418927638567)
		(3.4653061224489794,2.336580957877763)
		(3.561224489795918,2.3491532471578775)
		(3.657142857142857,2.3616529181824055)
		(3.753061224489796,2.3740752636690354)
		(3.8489795918367347,2.3864165264667445)
		(3.9448979591836735,2.3986737464491834)
		(4.040816326530612,2.410844633313496)
		(4.136734693877551,2.422927460516775)
		(4.2326530612244895,2.4349209765382342)
		(4.328571428571428,2.446824330403507)
		(4.424489795918367,2.458637008996773)
		(4.520408163265306,2.4703587841528787)
		(4.616326530612245,2.481989667893006)
		(4.7122448979591836,2.4935298744643695)
		(4.808163265306122,2.5049797880831117)
		(4.904081632653061,2.516339935472254)
		(5.0,2.52761096244274)
	}
	;
		\addplot[only marks, mark={square*}]
	coordinates {
			(1.3551020408163266,2.1016335501065315)
	}
	;
	\addlegendentry{$J(\alpha)$}
	\addlegendentry{$(\alpha_*, J(\alpha_*))$}
\end{axis}
\end{tikzpicture}
\caption{Plot of the function $\alpha \mapsto J(\alpha) := \frac{\alpha}{\gamma(\alpha)}$, for $b(x) = 
-x$ and $f(x) =x^2$. We prove in Proposition~\ref{prop:toy-model-multi-stab} that this function is 
decreasing on $(0, \alpha_*]$ and increasing on $[\alpha_*, \infty)$.  
}
\label{fig:plot toy model decreasing increasing}
\end{figure}
\begin{proposition}
	\label{prop:toy-model-multi-stab}
	Let $f(x) = x^2$ and $b(x) = -x$. There exists $\alpha_* > 0$ such that the function $\alpha 
	\mapsto 
	\frac{\alpha}{\gamma(\alpha)}$ is decreasing on $(0, \alpha_*]$ and increasing on 
	$[\alpha_*, \infty)$.  Moreover, one has
	\[ \lim_{\alpha \downarrow 0}{\frac{\alpha}{\gamma(\alpha)}} = +\infty, \quad \text{ and } \quad 
	\lim_{\alpha \rightarrow \infty}{\frac{\alpha}{\gamma(\alpha)}} = +\infty.  \]
	Let $J_* := \frac{\alpha_*}{\gamma(\alpha_*)}.$ We deduce that
	\begin{enumerate}
		\item For $J \in [0, J_*)$, $\delta_0$ is the unique invariant probability measure of \eqref{eq:McKeanVlasov}.
		\item For $J \in (J_*, \infty)$, \eqref{eq:McKeanVlasov} has three invariant probability measures: $\{\delta_0, 
		~\nu^\infty_{\alpha_1}, ~\nu^\infty_{\alpha_2}\}.$ with $\alpha_1 < \alpha_* < \alpha_2$.
		\item For $J = J_*$, \eqref{eq:McKeanVlasov} has two invariant probability measures: $\delta_0$ and 
		$\nu^\infty_{\alpha_*}$.
	\end{enumerate}
\end{proposition}
\begin{proof}
	The graph of the function $\alpha \mapsto \frac{\alpha}{\gamma(\alpha)}$ is plotted 
	Figure~\ref{fig:plot toy model decreasing increasing}. Define
	\begin{equation}
		\label{eq:definition de V} 
		\forall \alpha \geq 0,\quad V(\alpha) := \alpha \int_0^1{ (1+x) e^{\alpha^2 w(x)} \d x}. 
	\end{equation}
	\textbf{Claim} It holds that for all $\alpha > 0$
	\begin{equation}
	\label{eq:for b = -x, a formula for J}
	\frac{\alpha}{\gamma(\alpha)} = \frac{1}{\alpha} + V(\alpha).
	\end{equation}
	\textit{Proof}. 
	First, note that $\frac{\alpha}{\gamma(\alpha)} = \alpha \widehat{H_\alpha}(0)$. 
	In addition, by Lemma~\ref{lem:formula for H-K-theta-z} with $\theta = 0$ and $z = 0$, we find that $\widehat{H_\alpha}(0) = \frac{1}{\alpha^2} +  V(\alpha) / \alpha$. \qed

	Define for all $x \in [0,1)$
	\[ A(x) := \frac{-4 w(x)}{x^3} - (1+x) = \frac{1}{3} + 4 x^2 \sum_{k \geq 0} \frac{x^k}{k+5}. \]
	\textbf{Claim} It holds that
	\[ V'(\alpha) = \int_0^1{ A(x) e^{\alpha^2 w(x) } \d x}.  \]
	In particular $V$ is strictly increasing on $\mathbb{R}_+$. \\
	\textit{Proof}. 
	We have
	\begin{align*}
	V'(\alpha) = \int_0^1{ (1+x) e^{\alpha^2 w(x)} \d x} + 2 \alpha^2 \int_0^1{ (1+x) w(x) e^{\alpha^2 w(x) } 
	\d x}.
	\end{align*}
	Let
	\[ \forall x \in (0,1),\quad \theta(x) := \frac{(1+x)w(x)}{w'(x)}.\]
	We have  $\frac{w(x)}{w'(x)} = - \frac{(1-x) w(x) }{x^2}$ and so $\theta(x) = - \frac{1-x^2}{x^2} 
	w(x)$. In particular, $\theta$ can be extended to a $\mathcal{C}^1([0,1])$ function with $\theta(0) 
	= \theta(1) = 0$. Integrating par parts, we find that
	\begin{align*} 2 \alpha^2 \int_0^1{ (1+x) w(x) e^{\alpha^2 w(x) } \d x} &= 2 \alpha^2 \int_0^1{ 
	\theta(x) w'(x) e^{\alpha^2 w(x) } \d x} \\
	& = -2 \int_0^1{\theta'(x) e^{\alpha^2 w(x)} \d x}.
	\end{align*}
	Moreover, we have $\theta'(x) = \frac{2}{x^3} w(x) + (1+x)$ and so $ (1+x) - 2 \theta'(x) = A(x)$. \qed

	For all $\alpha \geq 1$, we have
	\[ V'(\alpha) \geq \frac{1}{3}  \int_0^1{  e^{\alpha^2 w(x) } dx} \geq \frac{1}{6 \alpha } \alpha  \int_0^1{ 
	(1+x) e^{\alpha^2 w(x) } dx} = \frac{1}{6 \alpha} V(\alpha). \]
	Consequently, we have $\forall \alpha \geq 1, ~V(\alpha) \geq V(1) \alpha^{1/6}$. Using \eqref{eq:for 
	b = -x, a formula for J}, we deduce that
	\[ \lim_{\alpha \downarrow 0}{\frac{\alpha}{\gamma(\alpha)}} = +\infty, \quad \text{ and } \quad 
	\lim_{\alpha \rightarrow \infty}{\frac{\alpha}{\gamma(\alpha)}} = +\infty.  \]
	It remains to study the variations of $\alpha \mapsto \frac{\alpha}{\gamma(\alpha)}$. Using 
	\eqref{eq:for b = -x, a formula for J}, we have
	\[ \frac{\d}{\d \alpha} \frac{\alpha}{\gamma(\alpha)} = \frac{\alpha^2 V'(\alpha) -1} {\alpha^2} = 
	\frac{W(\alpha^2) -1 }{\alpha^2}, \quad \text{ with } \quad  
	 W(\alpha) := \alpha \int_0^1{ A(x) e^{\alpha w(x)} \d x}. \]
	\textbf{Claim} The function $W$ is increasing on $\mathbb{R}_+$.\\
	\textit{Proof}. Let $D(x) := \frac{A(x) w(x)}{w'(x)}$. We have
	\begin{align*}
	W'(\alpha) &=  \int_0^1{ A(x) e^{\alpha w(x)} \d x} +  \int_0^1{ D(x) \alpha w'(x) e^{\alpha w(x)} \d x}  \\
	& =  \int_0^1{ [A(x) - D'(x)]  e^{\alpha w(x)} \d x}.
	\end{align*}
	To conclude it suffices to show that for all $x \in [0,1), ~ A(x) - D'(x) \geq 0$, which follows from 
	the explicit  formula satisfied by $A$ and $D$. \qed

	Finally, we have $\lim_{\alpha \rightarrow \infty} W(\alpha) = + \infty$.
	This follows from 
	\[ W(\alpha^2) = \alpha^2 V'(\alpha) \geq \alpha^2 \frac{1}{6 \alpha} V(1) 
	\alpha^{1/6}. \]
	Putting altogether, we deduce the result.
\end{proof}
\subsection{Conjecture on their stability}
Let $J \in (J_*, \infty)$. By Proposition~\ref{prop:toy-model-multi-stab}, \eqref{eq:McKeanVlasov} has exactly three invariant probability measures: $\{\delta_0, \nu^\infty_{\alpha_1}, \nu^\infty_{\alpha_2}\}$ with $\alpha_1 < \alpha_* < \alpha_2$. It is known that $\delta_0$ is attractive, see \cite{zbMATH06706511}. The question of the stability of $\nu^\infty_{\alpha_1}$ and $\nu^\infty_{\alpha_2}$ is more delicate. In view of Theorem~\ref{th:stability}, the stability is determined by the location of the zeros of 
\[ F(\alpha, z) := \widehat{H}_\alpha(z) - \widehat{\Psi}_\alpha(z). \] The explicit expression of $F(\alpha, z)$ is given in Proposition~\ref{prop:expression-of-Psi-and-H-hat} above. Recall the definition of $\lambda'_\alpha$ \eqref{eq:spectral-assumption}:
\[ -\lambda'_\alpha = \sup \{ \Re(z); \quad  F(\alpha, z) = 0 \}. \]
\begin{conjecture}
	We conjecture that $\lambda'_{\alpha_1} > 0$ and that $\lambda'_{\alpha_2} < 0$. \end{conjecture}
	In view of Theorem~\ref{th:stability}, this suggests that $\nu^\infty_{\alpha_1}$ is unstable and that $\nu^\infty_{\alpha_2}$ is stable. This conjecture is motivated by numerical investigations, see Figure~\ref{fig:fig}, and by the following analysis for $\alpha$ close to $\alpha_*$. First we note that for all $\alpha > 0$, it holds that:
\[ F(\alpha, 0) = \frac{\d}{\d \alpha} \frac{\alpha}{\gamma(\alpha)}. \]
In particular, for $\alpha = \alpha_*$, we have
\[ F(\alpha_*, 0) = 0. \]
The function $(\alpha, z) \mapsto F(\alpha, z)$ is $C^1$ in the neighborhood of $(\alpha_*, 0)$. In addition, we find that
\[ \partial_z F(\alpha_*, 0) > 0 \quad \text{ and } \quad \partial_\alpha F(\alpha_*, 0) > 0. \]
Therefore, the implicit function theorem applies, and gives the existence of a function $\alpha \mapsto z(\alpha)$ in the neighborhood of $\alpha_*$ such that $F(\alpha, z(\alpha)) = 0$. In addition, we have
\[ \frac{\d }{\d \alpha} z(\alpha_*) = - \frac{\partial_\alpha F(\alpha_*, 0)}{\partial_z F(\alpha_*, 0)} < 0. \]
This implies that for $\alpha < \alpha_*$, $\alpha$ sufficiently close to $\alpha_*$, it holds that $z(\alpha) > 0$ and so $\lambda'_\alpha < 0$.
Note however that when $\alpha < \alpha_*$, we have $z(\alpha) < 0$ but this local analysis is not sufficient to conclude that $\lambda'_\alpha > 0$. Indeed, there might be other solutions to the equation $F(\alpha, z) = 0$ with $\Re(z) \geq 0$ and $|z|$ far from zero.
\begin{figure}
\begin{subfigure}{0.45\textwidth} 
  \centering
  \begin{scaletikzpicturetowidth}{0.9\textwidth} 
\input{fig-1.107664180759389.tex}
  \end{scaletikzpicturetowidth}
  \caption{$F(\alpha_1, z) = 0$}
  \label{fig:sub-first}
\end{subfigure}
\begin{subfigure}{.45\textwidth}
  \centering
  \begin{scaletikzpicturetowidth}{0.9\textwidth}
\input{fig-1.7383261357223079.tex}
\end{scaletikzpicturetowidth}
\caption{$F(\alpha_2, z) = 0$}
  \label{fig:sub-second}
\end{subfigure}
\newline
\begin{subfigure}{.45\textwidth}
  \centering
  \begin{scaletikzpicturetowidth}{0.9\textwidth}
\input{fig-invm-1.107664180759389.tex}
\end{scaletikzpicturetowidth}
\caption{$\nu^\infty_{\alpha_1}(x)$}
  \label{fig:sub-third}
\end{subfigure}
\begin{subfigure}{.45\textwidth}
  \centering
  \begin{scaletikzpicturetowidth}{0.9\textwidth}
	\input{fig-invm-1.7383261357223079.tex}
\end{scaletikzpicturetowidth}
\caption{$\nu^\infty_{\alpha_2}(x)$}
  \label{fig:sub-fourth}
\end{subfigure}
\caption{
	Let $f(x) = x^2$ and $b(x) = -x$.
	For $J = 2.12$, the invariant probability measures of \eqref{eq:McKeanVlasov} are $\{\delta_0, \nu^\infty_{\alpha_1}, \nu^\infty_{\alpha_2} \}$ with $\alpha_1 \approx 1.108$ and $\alpha_2 \approx 1.7383$.
	The shape of the non-trivial invariant probability measures $\nu^\infty_{\alpha_1}$ and $\nu^\infty_{\alpha_2}$ are reported in Figures \ref{fig:sub-third} and \ref{fig:sub-fourth} respectively. 
	Figure~\ref{fig:sub-first}, we plot the curves $\Re F(\alpha_1, z) = 0$ (in blue) and $\Im F(\alpha_1, z)$ (in red). The two curves intersect at a zero of $F(\alpha_1, \cdot)$. We find numerically that $F(\alpha_1, 0.3065) \approx 0$. This suggests that $\nu^\infty_{\alpha_1}$ is unstable. For $\alpha = \alpha_2$, we find in Figure~\ref{fig:sub-second} that the zeros of $z \mapsto F(\alpha_2, z)$ have negative real part, suggesting that $\nu^\infty_{\alpha_2}$ is stable.
}
\label{fig:fig}
\end{figure}


\clearpage


\bibliographystyle{amsplainhyper_m}  
\bibliography{biblioquentincormier}

%


\ACKNO{The author thanks the referees for their useful comments. He also thanks Etienne Tanré, Romain Veltz and Eva Löcherbach for many valuable suggestions at different stages of this work.}

\nocite{*}
\end{document}